%% file: asynchronous.tex
\pdfoutput=1
\documentclass[12pt,draftcls,onecolumn]{IEEEtran}
\usepackage{fancyhdr}
\usepackage[margin=1in]{geometry}
\usepackage{amsmath}
\usepackage{graphicx}
\usepackage{amsfonts}
\usepackage{amssymb}
\usepackage{fixltx2e}
\usepackage{enumerate}
\usepackage{changepage}
\usepackage{mathtools}
\usepackage{color}
\usepackage{epstopdf}
\usepackage{empheq}
\usepackage{algorithm}
\usepackage{algpseudocode}

\usepackage{caption}
\usepackage{subcaption}

\usepackage{tikz}
\usetikzlibrary{matrix,positioning}

\DeclareMathOperator*{\argmin}{\arg\!\min}

\mathtoolsset{showonlyrefs=true}

\newcommand{\0}{\emptyset}
\newcommand{\R}{\mathbb{R}}
\newcommand{\N}{\mathbb{N}}

\newcommand{\Sp}{\mathbb{S}}

\newcommand{\dx}{\nabla_{x}}
\newcommand{\dxi}{\nabla_{x_i}}
\newcommand{\dxj}{\nabla_{x_j}}
\newcommand{\dmu}{\nabla_{\mu}}

\DeclarePairedDelimiter\floor{\lfloor}{\rfloor}
\newcommand{\grad}{\nabla}
\newcommand{\comp}{\Lambda}
\newcommand{\dualset}{M}
\newcommand{\fminx}{f^*}
\newcommand{\preg}{\alpha}
\newcommand{\dreg}{\beta}
\newcommand{\reg}{\kappa}
\newcommand{\compreg}{\Lambda_{\reg}}

\newcommand{\bmax}[1]{\left\|#1\right\|_{2,\infty}}
\newcommand{\btup}{(2, \infty)}
\newcommand{\pfunc}{L^t_{\reg}}

\newcommand{\Hi}{H^t_{\reg}}

\newcommand{\lve}{L_{\reg}}
\newcommand{\lxve}{\nabla_{x}\lve}
\newcommand{\lmve}{\nabla_{\mu}\lve}
\newcommand{\hmve}{\hat{\mu}_{\reg}}

\newcommand{\hxve}{\hat{x}_{\reg}}

\newcommand{\hzve}{\hat{z}_{\reg}}

\newcommand{\plip}{L_p}

\newcommand{\updatei}{\theta_i}

\newcommand{\Ni}{\mathcal{N}}

\newcommand{\Rij}{R^i_j}

\newcommand{\gradgm}{M_g}

\newcommand{\xdistmax}{L_{x}}
\newcommand{\xdiam}{D_{x}}
\newcommand{\xnorm}{M_x}
\newcommand{\gradgjm}{M_{g_j}}
\newcommand{\maxmu}{M_{\mu}}
\newcommand{\maxgradf}{M_{f}}
\newcommand{\mmax}{\hat{M}}
\newcommand{\regerror}{\epsilon}

\newcommand{\cloudxt}{x^c_t}

\newcommand{\new}[1]{#1}

\newcommand{\dkt}{D(k_t)}
\newcommand{\qp}{q_p}
\newcommand{\qd}{q_d}

\newtheorem{definition}{Definition}
\newtheorem{example}{Example}
\newtheorem{lemma}{Lemma}

\newtheorem{proposition}{Proposition}
\newtheorem{theorem}{Theorem}
\newtheorem{talgorithm}{Algorithm}
\newtheorem{assumption}{Assumption}
\newtheorem{remark}{Remark}
\newtheorem{problem}{Problem}

\author{Matthew T. Hale$^\star$, Angelia Nedi\'{c}$^\dagger$, and Magnus Egerstedt$^\star$\thanks{$^\star$School of Electrical and Computer Engineering, Georgia Institute of
Technology, Atlanta, GA 30332, USA. Email: \texttt{\{matthale, magnus\}@gatech.edu}. Research supported in
part by  the NSF under Grant CNS-1239225.}\thanks{$^\dagger$School of Electrical, Computer and Energy Engineering, Arizona State University, Tempe, AZ 85287, USA. Email: \texttt{
Angelia.Nedich@asu.edu}.}
}

\begin{document}
\title{Asynchronous Multi-Agent Primal-Dual Optimization}
\maketitle

\input{0-abstract.tex}
\input{1-introduction.tex}

\input{2-preliminaries.tex}

\input{3-error.tex}

\input{4-decentralized.tex}

\input{5-convergence.tex}

\input{6-counterexample.tex}

\input{7-simulation.tex}

\input{8-conclusion.tex}

\bibliographystyle{IEEEtran}{}
\bibliography{sources}

\end{document}

%% file: 0-abstract.tex
\begin{abstract}
We present a framework for asynchronously solving convex optimization problems over networks of agents
which are augmented by the presence of a centralized cloud computer. 
This framework uses a Tikhonov-regularized primal-dual approach
in which the agents update the system's primal variables and the cloud updates
its dual variables. 
To minimize coordination requirements placed upon the system,
the times of
communications and computations among the agents are allowed to be arbitrary, provided
they satisfy mild conditions. Communications
from the agents to the cloud are likewise carried out without any coordination in their timing. 
However, we require that the cloud keep the dual variable's value synchronized across
the agents, and a counterexample is provided that demonstrates that 
this level of synchrony is indeed necessary for convergence. Convergence rate estimates are provided
in both the primal and dual spaces, and simulation results are presented that demonstrate
the operation and convergence of the proposed algorithm. 
\end{abstract}

%% file: 1-introduction.tex
\section{Introduction}
Networked coordination and optimization have been
applied across a broad range of 
application domains, such as
sensor networks \cite{khan09,trigoni12,cortes02,rabbat04}, robotics \cite{soltero13}, 
smart power grids \cite{vytelingum10,caron10}, and
communications \cite{kelly98,chiang07,mitra94}.
A common feature of some applications is
the (sometimes implicit) assumption that communications and computations
occur in a synchronous fashion. More precisely, though no agent may have access to all information
in a network, the information it does have access to
is assumed to be up-to-date and/or computations onboard the agents
are assumed to occur concurrently. 

One can envision several reasons why these synchrony assumptions may fail.
Communications may interfere with each other, slowing
data transmissions, or else they may occur serially over a shared
channel, resulting in delays when many messages must be sent.
In other cases it may simply be undesirable
to stay in constant communication in a network due to the energy
required to do so, e.g., within a team of battery-powered robots. 
Apart from communication delays, it may be the case that some agents
produce new data, such as a new state value, faster than other agents,
leading to mismatches in update rates and thus
mismatches in when information becomes available. 
Regardless of their cause, 
the resulting delays are often 
unpredictable in duration, and
the timeliness of any piece of information in a network 
with such delays typically cannot
be guaranteed. 
While one could simply have agents pause their computations
while synchronizing information across a network, it has been
shown that asynchronous algorithms \new{can} outperform
their synchronous counterparts which pause to synchronize information \cite[Section 6.3.5]{bertsekas97}\cite[Section 3.3]{bertsekas89}.
Accordingly, this paper focuses on asynchronous algorithms for multi-agent optimization.

In particular, this paper considers
multi-agent convex optimization problems that need not be separable, 
and its structural novelty comes from the introduction of a centralized cloud computer
\new{and its associated communications model}. The cloud's role is to aggregate centralized information
and perform centralized computations for the agents in the network, and
the motivation for including a cloud computer comes from its 
ability to communicate with many devices and its ability to 
provide ample processing power remotely. 
However, the cloud's operations take time to perform
specifically because they are centralized, 
and although the cloud adds
centralized information to a network, 
the price one has to pay for this centralized information
is that it is generated slowly.
As such, the proposed algorithmic model has to take this slowness into account.

In this paper we consider problems in which each agent has a local cost
and local set constraint, and in which the network itself is associated with a 
 non-separable
coupling cost. 
The agents are moreover subject to non-separable ensemble-level
inequality constraints\footnote{The work here can include equality constraints without
any further changes, though we focus only on inequality constraints for notational simplicity.} that could, for example, correspond to shared resources.
To solve these problems, we consider a primal-dual approach that allows the agents' behavior to be totally asynchronous \cite[Chapter 6]{bertsekas97} 
(cf. partially asynchronous \cite[Chapter 7]{bertsekas97}),
both when communicating among themselves and when transmitting to the cloud. 
However, we do require that the cloud's transmissions to the agents always keep the dual
variable's value synchronized among the agents. This synchrony is verified to be necessary
in Section~\ref{sec:counter}, where
  a counterexample shows that allowing the agents
to disagree upon the value of the system's dual variable 
can preclude convergence altogether.
The dual variable's value is the lone point of synchrony in the presented algorithm, 
and all other aspects of the system are designed to strive toward operating as asynchronously as possible
in a general optimization setting. 

To produce such an algorithm, we apply a Tikhonov regularization
to the Lagrangian associated with the problem of interest.
\new{
This regularization causes the algorithm to only approximately solve optimization problems, and error bounds
are provided in terms of the regularization parameters, along with a choice rule for selecting
these parameters to enforce any desired error bound. The regularization we use induces a tradeoff between
speed and accuracy in the optimization process, and it is shown that requiring a less accurate
solution allows the algorithm to converge faster and vice versa. 

}We also make use of an existing framework for asynchronous optimization \cite[Sections 6.1-6.2]{bertsekas97}\cite{bertsekas83},
which accommodates general unconstrained or set-constrained problems. This framework hinges upon
the ability to construct a sequence of sets satisfying certain properties which
admit a Lyapunov-like convergence result, and we show that our regularization guarantees the ability to construct
this sequence of sets as long as the problem satisfies mild assumptions. 
We also provide novel convergence rate estimates in both the primal and dual spaces that explicitly account
for the delays in the system. The contribution of this work thus consists of an
asynchronous primal-dual optimization algorithm together with its convergence rates. 

There exists a large corpus of work on multi-agent optimization that is related
to the work here. In \cite{bertsekas97} a range of results are gathered on
asynchronous multi-agent optimization (for problems without functional constraints
or with linear equality constraints) 
in Chapters 6 and 7. Earlier work on asynchronous algorithms can be traced back to \cite{chazan69} and \cite{baudet78},
which consider fixed points of certain classes of operators. 
Long-standing optimization algorithms known as the
Jacobi and Gauss-Seidel methods are also covered in \cite{bertsekas97} for linear problems in Section 2.4.
Linear consensus type problems are studied in \cite{olfati-saber04},
including cases in which identical time delays are associated with the communication channels. 
The framework in \cite[Sections 6.1-6.2]{bertsekas97} is the most general, and we therefore
use it as our starting point for optimization in the primal space. 

\new{A key difference between
our work and earlier work is that we asynchronously solve
general constrained convex optimization problems which, in general, need not
satisfy the conditions in \cite{bertsekas97,chazan69,baudet78}. The work in \cite{zhu12} also solves
constrained optimization problems asynchronously, though it requires bounded communication delays
between agents and has each agent updating both a full primal vector and a full dual vector. In the
current paper, communication delays do not have a uniform bound, and each agent updates only its
own state as would be the case, e.g., in a team of robots. 
}

Two other relevant and well-known algorithms of current interest 
are gossip algorithms and the alternating direction method of multipliers (ADMM).
Here we do not consider gossip-type algorithms since they either require synchronous communications among
the agents or, in the asynchronous case, allow only one communication channel to be active at a time \cite{boyd06},
and our aim is to support communication models that are as
general as possible \new{by allowing any number of links to be active at a time}.

In contrast to this, ADMM essentially imposes a Gauss-Seidel structure among the primal updates 
made by the agents \cite{boyd11}. 
Related work in \cite{zhang14} presents an asynchronous variant of ADMM, though 
it requires bounded delays and updates of all primal and dual variables onboard each agent,
neither of which are required here. 
The algorithm we present can be viewed
as a method related to ADMM that allows all agent behaviors to be essentially arbitrary
in their timing. This provides a great degree of flexibility in the agents' primal updates by
not requiring any particular ensemble update rule or bounded delays, or requiring an agent to 
update all variables in the system. 

The remainder of the paper is organized as follows. Section~\ref{sec:problem} describes
the optimization problem to be solved and the regularization used. Then
Section~\ref{sec:error} gives a rule for choosing regularization parameters to limit
errors in the system. Next, Section~\ref{sec:decent} provides the asynchronous
algorithm that is the main focus of the paper. Then
Section~\ref{sec:main} proves convergence of the asynchronous algorithm 
and provides convergence rates for it. We show in Section~\ref{sec:counter} that synchrony
in the dual variable is indeed a necessary condition for convergence.
Section~\ref{sec:simulation} presents simulation results for the asynchronous algorithm, and
 Section~\ref{sec:con} concludes the paper.

%% file: 2-preliminaries.tex
\section{Multi-Agent Optimization} \label{sec:problem}
This section gives a description of the problems under consideration 
and establishes key notation. 
To that end, the symbol $\|\cdot\|$ 
without a subscript always denotes the Euclidean norm. 
We use the notation $x_{-i}$ to denote the vector $x$ with its
$i^{th}$ component removed, i.e.,
\begin{equation}
{x_{-i} = (x_1, \ldots, x_{i-1}, x_{i+1}, \ldots, x_n)}. 
\end{equation} 
We also define the index set $[P] := \{1, \ldots, P\}$
for all $P \in \N$, and we will use the term ``ensemble'' to refer
to aspects of the problem that involve all agents. 

\subsection{Problem Statement} \label{ss:problem}
This paper solves convex optimization problems
over networks comprised by $N$ agents. The agents are indexed over $i \in [N]$,
and agent $i$ has an associated decision variable, $x_i \in \R^{n_i}$, with $n_i \in \N$, and we allow
for $n_i \neq n_j$ when $i \neq j$. Each agent has to satisfy a local set constraint,
expressed by requiring
$x_i \in X_i \subset \R^{n_i}$,
where we assume the following about each $X_i$.

\begin{assumption} \label{as:sets}
For all $i \in [N]$, the set $X_i$ is non-empty, compact, and convex. \hfill $\lozenge$
\end{assumption}

Note that Assumption~\ref{as:sets} allows for box constraints, which are common in 
multi-agent optimization.
We will also refer to the ensemble decision variable of the network, defined as
$
{x = (x_1^T,\ldots,x_N^T)^T \in X := X_1 \times \cdots \times X_N \subset \R^n},
$
where $n = \sum_{i \in [N]} n_i$. 
Assumption~\ref{as:sets} guarantees that $X$ is also non-empty, compact, and convex. 

Agent $i$ seeks to minimize a local objective function $f_i : X_i \to \R$ which depends
only upon $x_i$. 
Together, the agents also seek to minimize a coupling cost $c : \R^n \to \R$ which depends
upon all states and can be non-separable. 
We impose the following assumption on $c$ and each $f_i$. 
\begin{assumption} \label{as:objectives}
For all $i \in [N]$, the function $f_i$ is convex and $C^2$ \new{(twice continuously differentiable)} in $x_i$. 
The function $c$ is convex and $C^2$ in $x$. 
\hfill $\lozenge$
\end{assumption}

Gathering these costs gives 
\begin{equation}
f(x) = c(x) + \sum_{i \in [N]} f_i(x_i), 
\end{equation}
and when Assumptions~\ref{as:sets} and \ref{as:objectives} hold, $f$ has a 
well-defined minimum value over $X$. 

We consider problems with
ensemble-level inequality constraints, namely we require that the inequality
\begin{equation}
g(x) := \big(g_1(x), g_2(x), \ldots, g_m(x)\big)^T \leq 0
\end{equation}
hold component-wise, under the following assumption.
\begin{assumption} \label{as:constraints}
The function $g : \R^n \to \R^m$ is convex and $C^2$ in $x$. \hfill $\lozenge$
\end{assumption}

In particular, $g$ does not need to be separable. 
At the ensemble level, 
we now have a convex optimization problem, stated below.
\begin{problem} \label{prob:ensemble}
\begin{align}
\textnormal{minimize } &f(x) \\
\textnormal {subject to } &g(x) \leq 0\\
                      &x \in X. \tag*{$\blacklozenge$}
                      \end{align} 
\end{problem}

On top of the problem formulation itself, Section~\ref{sec:decent} will 
specify an architecture
that provides a mixture of distributed information sharing among 
agents and centralized information from a cloud computer. 
As a result,
the solution to Problem~\ref{prob:ensemble} will involve 
a distributed primal-dual algorithm since such an algorithm 
 is implementable in a natural way on the cloud-based architecture.
Towards enabling this algorithm, 
 we enforce Slater's condition \cite[Assumption 6.4.2]{bertsekas03} 
which enables us to find a compact set that contains the optimal dual point in
Problem~\ref{prob:ensemble}.
\begin{assumption} (\emph{Slater's condition}) \label{as:slater}
There exists a point $\bar{x} \in X$ such that $g(\bar{x}) < 0$. \hfill $\lozenge$
\end{assumption}

\subsection{An Ensemble Variational Inequality Formulation} \label{ss:ensemble}
Under Assumptions~\ref{as:sets}-\ref{as:slater} we define an ensemble 
variational inequality 
in terms of Problem~\ref{prob:ensemble}'s Lagrangian. 
The Lagrangian associated with Problem~\ref{prob:ensemble} is defined as
\begin{equation}
L(x, \mu) = f(x) + \mu^Tg(x), 
\end{equation}
where $\mu \in \R^m_{+}$ and $\R^m_{+}$ denotes the non-negative orthant of $\R^m$. 
By definition, $L(\cdot, \mu)$ is convex for all $\mu \in \R^m_{+}$
and $L(x, \cdot)$ is concave for all $x \in X$. These properties and the differentiability assumptions placed upon $f$ and $g$ together
imply that $\nabla_x L(\cdot, \mu):= \frac{\partial L}{\partial x}(\cdot, \mu)$ and $-\nabla_{\mu} L(x, \cdot):= -\frac{\partial L}{\partial \mu}(x, \cdot)$ 
are monotone operators on their respective domains.
It is known that Assumptions~\ref{as:sets}-\ref{as:slater} 
imply that a point $(\hat{x}, \hat{\mu}) \in X \times \R^m_{+}$ is a solution to Problem~\ref{prob:ensemble}
if and only if it is a saddle point of $L$ \cite{kuhn51}, i.e., it maximizes $L$ over $\mu$ and minimizes $L$ over $x$ so that it satisfies the inequalities
\begin{equation} \label{eq:saddle}
L(\hat{x}, \mu) \leq L(\hat{x}, \hat{\mu}) \leq L(x, \hat{\mu})
\end{equation}
for all $x \in X$ and $\mu \in \R^m_{+}$. From Assumptions~\ref{as:sets}-\ref{as:slater} it is guaranteed that a saddle point $(\hat{x}, \hat{\mu})$ exists
\cite[Corollary 2.2.10]{facchinei03}. 

Defining the symbol $\hat{z}$ to denote a saddle point via
$\hat{z} = (\hat{x}, \hat{\mu})$ and using 
$z = (x, \mu)$ to denote an arbitrary point in $X \times \R^m_{+}$,
we define the composite gradient operator
\begin{equation}
\comp(z) := \comp(x, \mu) = \left(\begin{array}{r} \nabla_x L(x, \mu) \\ -\nabla_{\mu} L(x, \mu) \end{array}\right).
\end{equation}
Then
the saddle point condition in Equation \eqref{eq:saddle} can be restated as the following
ensemble variational inequality \cite[Page 21]{facchinei03}.
\begin{problem} \label{prob:ensemblevi}
Find a point $\hat{z} \in X \times \R^m_{+}$ such that
${(z - \hat{z})^T\comp(\hat{z}) \geq 0}$
for all $z \in X \times \R^m_{+}$. \hfill $\blacklozenge$
\end{problem}

\subsection{Tikhonov Regularization} \label{ss:reg}
Instead of solving Problem~\ref{prob:ensemblevi} as stated, 
we regularize the problem in order to make it more readily solved
asynchronously and to enable us to analyze the convergence rate of
the forthcoming asynchronous algorithm. 
The remainder of this paper will make extensive use of this regularization.
First we define $\eta$-strong convexity for a differentiable function. 

\begin{definition} \label{def:strongcon}
A differentiable function $f$ is said to be $\eta$-strongly convex if
\begin{equation}
\big(\nabla f(v_1) - \nabla f(v_2)\big)^T(v_1 - v_2) \geq \eta\|v_1 - v_2\|^2
\end{equation}
for all $v_1$ and $v_2$ in the domain of $f$. \hfill $\Diamond$
\end{definition}

We now regularize the Lagrangian 
using constants $\preg > 0$ and $\dreg > 0$ to get
\begin{equation}
L_{\preg, \dreg}(x, \mu) = f(x) + \frac{\preg}{2}\|x\|^2 + \mu^Tg(x) - \frac{\dreg}{2}\|\mu\|^2,
\end{equation}
where we see that $L_{\preg, \dreg}(\cdot, \mu)$ is $\preg$-strongly convex and $L_{\preg, \dreg}(x, \cdot)$ is $\dreg$-strongly concave (which is equivalent
to $-L_{\preg, \dreg}(x, \cdot)$ being $\dreg$-strongly convex).
Accordingly, $\dx L_{\preg, \dreg}(\cdot, \mu)$ and $-\dmu L_{\preg, \dreg}(x, \cdot)$ are strongly monotone operators over their domains. 
We also define $\reg = (\preg, \dreg)$ and replace the subscripts $\preg$ and $\dreg$ with the single subscript
$\reg$ for brevity when we are not using specific values of $\alpha$ and $\beta$. 
We now have the regularized composite gradient operator,
\begin{equation} 
\compreg(z) := \compreg(x, \mu) = \left(\begin{array}{r} \lxve(x, \mu) \\ -\lmve(x, \mu) \end{array}\right) : X  \times \R^m_{+} \to \R^{m + n}.
\end{equation}
The strong monotonicity of $\nabla_x \lve(\cdot, \mu)$ and $-\nabla_{\mu} \lve(x, \cdot)$ together imply
that $\compreg$ itself is strongly monotone, and Assumptions~\ref{as:sets}-\ref{as:slater} imply
that $\lve$ has a unique saddle point, $\hzve$ \cite[Theorem 2.3.3]{facchinei03}. 
We now focus on solving the following regularized ensemble variational inequality.
\begin{problem} \label{prob:regvi}
Find the point $\hzve := (\hxve, \hmve) \in X \times \R^m_{+}$ such that
${(z - \hzve)^T\compreg(\hzve) \geq 0}$
for all $z \in X \times \R^m_{+}$. \hfill $\blacklozenge$
\end{problem}

As a result of the regularization of $\comp$ to define $\compreg$,
the solution $\hzve$ will not equal $\hat{z}$. In particular, for a solution $\hat{z} = (\hat{x}, \hat{\mu})$ to 
Problem~\ref{prob:ensemblevi} and the solution $\hzve = (\hxve, \hmve)$ to Problem~\ref{prob:regvi}, we will have
$\hat{x} \neq \hxve$ and $\hat{\mu} \neq \hmve$. 
Thus the regularization done with $\preg$ and $\dreg$ affords us a greater ability to find
saddle points asynchronously and, as will be shown, the ability to estimate convergence rates towards a solution, but does so
at the expense of accuracy by changing the solution itself. While solving Problem~\ref{prob:regvi}
does not result in a solution to Problem~\ref{prob:ensemblevi}, the continuity of $\lve$ over $X \times \R^m_{+}$ suggests
that using small values of $\preg$ and $\dreg$ should lead to small differences between $\hat{z}$ and $\hzve$
so that the level of error introduced by regularizing is acceptable in many settings. 
Along these lines, we provide a choice rule for $\preg$ and $\dreg$ in Section~\ref{sec:error} that enforces any
desired error bound for certain errors due to regularization. 

There is a well-established literature regarding projection-based methods for solving variational inequalities like 
that in Problem~\ref{prob:regvi}, e.g., \cite[Chapter 12.1]{facchinei03}.
We seek to use projection methods because they naturally fit with the mixed centralized/decentralized architecture to be covered 
in Section~\ref{sec:decent}, though 
it is required that $\compreg$ be Lipschitz to make use of such methods. Currently, 
$\compreg$ cannot be shown to be Lipschitz
because its domain, $X \times \R^m_{+}$, is unbounded. To rectify this situation, we now determine a non-empty, compact, convex set $\dualset \subset \R^m_{+}$ 
which contains $\hmve$, allowing us to solve Problem~\ref{prob:regvi} over a compact domain. 
Below, we use the unconstrained minimum value of $f$ over $X$,
$
\fminx  := \min_{x \in X} f(x),
$
which is well-defined under Assumptions~\ref{as:sets} and \ref{as:objectives}. 
We have the following result based upon \cite[Chapter 10]{uzawa58a}. 

\begin{lemma} \label{lem:dualset}
Let $\bar{x} \in X$ be a Slater point of $g$. Then 
\begin{equation}
\hmve \in \dualset := \left\{\mu \in \R^m_{+} : \|\mu\|_1 \leq \frac{f(\bar{x}) + \frac{\preg}{2}\|\bar{x}\|^2 - \fminx}{\min\limits_{1 \leq j \leq m} \{-g_j(\bar{x})\}}\right\}.
\end{equation}
\end{lemma}
\emph{Proof:} See \cite{hale15}, Section II-C. \hfill $\blacksquare$

\new{
If $\fminx$ is not available, any lower bound on $\fminx$ can be used in defining
$\dualset$, and the above construction is still valid when using such a lower bound 
in conjunction with any Slater point $\bar{x} \in X$.}
Having defined $M$, we see that the norm of the gradient of $\dx \lve(\cdot, \mu)$ 
can be uniformly upper-bounded for all $\mu \in M$,
and $\dx \lve(\cdot, \mu)$ is therefore Lipschitz. Denote
its Lipschitz constant by $\plip$. 
We now define a synchronous, ensemble-level primal-dual projection method for
finding $\hzve$ based on \cite{uzawa58a}. It relies on the Euclidean projections onto
$X$ and $\dualset$, denoted $\Pi_X[\cdot]$ and $\Pi_{\dualset}[\cdot]$, respectively. 

\begin{talgorithm} \label{alg:ensemble}
Let $x(0) \in X$ and $\mu(0) \in \dualset$ be given. For values $k = 0,1,\ldots$, execute
\begin{align}
x(k+1)   &= \Pi_X\left[x(k) - \gamma\left(\nabla_{x}\lve\left(x(k),\mu(k)\right)\right)\right] \\
\mu(k+1) &= \Pi_{\dualset}\left[\mu(k) + \rho\left(\nabla_{\mu}\lve\left(x(k), \mu(k)\right)\right)\right]. \tag*{$\triangle$} 
\end{align}
\end{talgorithm}

Here $\gamma$ and $\rho$ are stepsizes whose values will be determined in Theorems~\ref{thm:dualconv} and \ref{thm:pconv} in Section~\ref{sec:main}. 
Algorithm~\ref{alg:ensemble} will serve as a basis for the asynchronous algorithm developed in Section~\ref{sec:decent}, though, 
as we will see, significant modifications must be made to this update law to account for asynchronous behavior
in the network.

%% file: 3-error.tex
\section{Bounds on Regularization Error} \label{sec:error} 
In this section we briefly cover bounds on two 
errors that result from the Tikhonov regularization of $L$. 
For more discussion, we refer the reader
to Section 3.2 in \cite{koshal11} for regularized Lagrangian
methods and to Chapter 12.2
in \cite{facchinei03} for a discussion of regularization
error in variational inequalities. 

For any fixed choice of $\preg$ and $\dreg$, denote the corresponding solution to Problem~\ref{prob:regvi} by $\hat{z}_{\preg, \dreg}$. 
It is known that as $\preg \downarrow 0$ and $\dreg \downarrow 0$ across a sequence of problems, the solutions $\hat{z}_{\preg, \dreg} \to \hat{z}_0$, 
where $\hat{z}_0$
is the solution to Problem~\ref{prob:ensemblevi} with least Euclidean norm \cite[Theorem 12.2.3]{facchinei03}.
Here, we are not interested in solving a sequence of problems for evolving values of $\preg$ and $\dreg$
because of the computational burden of doing so;
instead, we solve only a single problem. It is also known that an algorithm with an iterative regularization
wherein $\preg$ and $\dreg$ tend to zero as a function of the iteration number 
can also converge to $\hat{z}_0$ \cite{bakushinskii74}, 
though here it would be difficult to synchronize changes in 
the regularization parameters across the network. As a result, we proceed with a fixed regularization and give
error bounds in terms of the regularization parameters we use. 

Our focus is on selecting the parameters
$\preg$ and $\dreg$ to satisfy desired bounds on errors introduced
by the regularization. First we present error bounds and
then we cover how to select $\preg$ and $\dreg$ to bound
these errors by any positive constant. 

\subsection{Error Bounds}
Below we use the following four constants:
\begin{equation}
\maxgradf := \max_{x \in X} \|\grad f(x)\| \qquad \maxmu := \max_{\mu \in M} \|\mu\| \qquad \gradgjm := \max_{x \in X} \|\grad g_j(x)\| \qquad \xnorm := \max_{x \in X} \|x\|.
\end{equation}
We first state the error in optimal cost. 

\begin{lemma} \label{lem:costdiff}
Let Assumptions~\ref{as:sets}-\ref{as:slater} hold. 
For regularization parameters $\preg > 0$ and $\dreg > 0$, the error in optimal cost incurred 
by regularizing $L$ is bounded according to 
\begin{equation}
|f(\hat{x}_{\preg, \dreg}) - f(\hat{x})| \leq \maxgradf\maxmu\sqrt{\dreg/2\preg} + \frac{\preg}{2}\xnorm^2.
\end{equation}
\end{lemma}
\emph{Proof:} See \cite[Lemma 3.3]{koshal11}. \hfill $\blacksquare$

Next we bound the constraint violation that is possible in solving Problem~\ref{prob:regvi}. 

\begin{lemma} \label{lem:gviolation}
Let Assumptions~\ref{as:sets}-\ref{as:slater} hold. For $\preg > 0$ and $\dreg > 0$, the constraint violation
due to regularizing $L$ is bounded according to
\begin{equation}
\max\{0, g_j(\hat{x}_{\preg, \dreg})\} \leq \gradgjm \maxmu \sqrt{\dreg/2\preg} \textnormal{ for all } j \in [m].
\end{equation}
\end{lemma}
\emph{Proof:} See \cite[Lemma 3.3]{koshal11}. \hfill $\blacksquare$

\subsection{Selecting Regularization Parameters} 
We now discuss \new{one possible choice} rule for selecting $\preg$ and $\dreg$ based upon Lemmas~\ref{lem:costdiff} and \ref{lem:gviolation}. 
Both lemmas suggest using ${\dreg < \preg}$ to achieve smaller errors and, given that we expect $\preg < 1$,
we choose $\dreg = \preg^3/2$. Suppose that there is some maximum error $\regerror > 0$ specified for Lemmas~\ref{lem:costdiff}
and \ref{lem:gviolation}. The following result \new{provides sufficient conditions for enforcing
this bound} by choosing $\preg$ and $\dreg$ appropriately.

\begin{proposition} \label{prop:regchoose}
Let $\regerror > 0$ be given. For 
\begin{equation}
\mmax = \max\left\{\max_{j \in [m]} \gradgjm\maxmu, \maxgradf\maxmu\right\}, 
\end{equation}
choosing regularization parameters
${\preg < 2\regerror/(\mmax + \xnorm^2)}$ and ${\dreg = \preg^3/2}$
gives
\begin{equation}
\max\{0, g_j(\hat{x}_{\preg, \dreg})\} < \regerror \,\, \textnormal{ and } \,\, |f(\hat{x}_{\preg, \dreg}) - f(\hat{x})| < \regerror.
\end{equation}
\end{proposition}

\emph{Proof:} 
By definition of $\mmax$ and Lemmas~\ref{lem:costdiff} and \ref{lem:gviolation}, 
\begin{equation}
\max\{0, g_j(\hat{x}_{\preg, \dreg})\} \leq \mmax \sqrt{\dreg/2\preg} + \xnorm^2\preg/2
\end{equation}
for all $j \in [m]$ and
\begin{equation}
|f(\hat{x}_{\preg, \dreg}) - f(\hat{x})| \leq \mmax \sqrt{\dreg/2\preg} + \xnorm^2\preg/2.
\end{equation}
Thus we require
$\mmax \sqrt{\dreg/2\preg} + \xnorm^2\preg/2 < \regerror$. 
Choosing $\dreg = \preg^3/2$ 
and solving for $\preg$ gives the desired bound. 
\hfill $\blacksquare$

%% file: 4-decentralized.tex
\section{Asynchronous Optimization} \label{sec:decent}
In this section, we examine what happens when primal and dual updates
are computed asynchronously. The agents compute primal updates and
the cloud computes dual updates and, because the cloud is centralized, the dual
updates in the system are computed slower than the primal updates are.
Due to the difference in primal and dual
update rates, we will now index
the dual variable, $\mu$, over the time index $t$ and will continue
to index the primal variable, $x$, over the time index $k$. 
In this section we make use of the optimization framework in 
Sections 6.1 and 6.2 of \cite{bertsekas97}. 
Throughout this section,
discussions will have the same value of $\mu(t)$ onboard all agents
simultaneously, and this is shown to be a necessary condition for convergence
in Section~\ref{sec:counter}. 

\subsection{Per-Agent Primal Update Law}
The exact update law used by agent $i$ will be detailed below. For the
present discussion, we need only to understand a few basic facts
about the distribution of communications and computations in the system. 
Agent $i$ will store values of \new{some other agents' states} in its onboard
computer, but will only update its own state within that state vector;
states stored by agent $i$ corresponding to other agents 
will be updated only when those agents send their state values to agent $i$. 
Because these operations occur asynchronously, there is no 
reason to expect that agents $i$ and $j$ (with $i \neq j$) will agree
upon the values of any states in the network. 

As a result, we index
each agent's state vector using a superscript: agent $i$'s copy
of the state of the system is denoted $x^i$
and agent $i$'s copy of its own state is denoted $x^i_i$. 
In this notation we say that agent $i$ updates $x^i_i$ but not 
$x^i_j$ for any $j \neq i$. The state value $x^i_j$ is precisely the
content of messages from agent $j$ to agent $i$ and its value onboard
agent $i$ is changed only when agent $i$ receives messages from agent $j$
(and this change occurs immediately when messages are received by agent $i$). 

\new{
To prevent unnecessary communications among the agents,
we only require two agents to communicate if each needs the other's state value
in its computations. We make this notion precise in the following definition.

\begin{definition} \label{def:essential}
Agent $j$ is an \emph{essential neighbor} of agent $i$ (where $i \neq j$) if
$\dxj \lve := \frac{\partial \lve}{\partial x_j}$ depends upon $x_i$. 
The set of indices of all
essential neighbors of agent $i$ is called its \emph{essential neighborhood},
denoted $\Ni_i$. \hfill $\Diamond$
\end{definition}

We illustrate the role of Definition~\ref{def:essential} in defining
communications among the agents in the following example. 

\begin{example} \label{ex:essential}
Consider a system with four agents with scalar states. For all $i \in [4]$, we have $f_i(x_i) = x_i$,
and the constraints are $g_1(x) = \frac{1}{2}(x_1 - x_2)^2$ and $g_2(x) = \frac{1}{2}(x_3 - x_4)^2$, 
with $c \equiv 0$. 
For $\preg, \dreg > 0$,
we find that $\grad_{x_1} \lve(x, \mu) = (1 + \preg + \mu_1)x_1 - \mu_1x_2$. As a result, agent $1$'s
essential neighborhood is $\Ni_1 = \{2\}$. We also find $\Ni_2 = \{1\}$, 
$\Ni_3 = \{4\}$ and $\Ni_4 = \{3\}$. As a result, agents $1$ and $2$ need only to 
communicate with each other and store each other's states; neither needs to communicate with agents
$3$ or $4$ at any point, nor to store the states of agents $3$ and $4$. 
Similarly, agents $3$ and $4$ communicate and store each other's states, but never 
communicate with agents $1$ and $2$ and therefore do not store their states.
Agents that communicate states with each other do not need to do so simultaneously and 
can do so with any timing. 
Then at each timestep, there are four possible directed edges that can be active, and
the union of all communication graphs over all timesteps is 
shown in Figure~\ref{fig:exgraphs}. 

\begin{figure}
\centering
\includegraphics[width=1.5in]{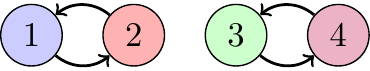}
\caption{
The union of all possible communication graphs
over all timesteps in Example 1. This graph is
neither complete, nor is it even (strongly or weakly)
connected. 
}
\label{fig:exgraphs}
\end{figure}

Here we see that the agents' communications need not
comprise a graph which is complete, nor even one which is connected in any sense. 
What results then is a system in which there may be multiple groups of agents which do not
interact at all and which may indeed not even know of each other's existence, though they are
jointly solving an optimization problem. \hfill $\triangle$
\end{example}

Clearly $j \in \Ni_i$ if and only if $i \in \Ni_j$, and thus agent $i$ both sends
information to and receives information from its essential neighbors. 
While each agent only needs to store the states of its essential neighbors, 
we proceed as though each agent stores a full state vector in order to
circumvent the need to track different dimensions of agents' states.
For agent $i$, one can assume that $x^i_j$ is fixed at zero for all
$j \not\in \Ni_i$. 
Rather than considering a fixed communication topology and analyzing an optimization algorithm
developed over that topology, Example~\ref{ex:essential} shows that we take the opposite approach:
the information dependencies in the system determine which agents must communicate because
these dependencies define the agents' essential neighborhoods. Of course, in some cases,
it will be difficult for two agents to communicate and they will do so only occasionally
and without any specified schedule, and this is permitted by the asynchronous
problem formulation we develop below. 
}

The agents' primal updates also do not occur concurrently with dual
updates (which will be computed by the cloud). It is therefore necessary to track
which dual variable the agents currently have onboard. 
At time $k$, if agent $i$ has $\mu(t)$ in its onboard computer, we denote
agent $i$'s copy of $x$ by $x^i(k; t)$.

Each agent is allowed to compute its state updates using any clock
it wishes, regardless of the timing of the other agents' clocks. We use
the symbol $K$ to denote a virtual global clock which contains
the clock ticks of each agent's clock, and $K$ can be understood as 
containing ordered indices of instants in time at which some number of agents compute
state updates.
Without loss of generality we take $K = \N$. 
We denote the set of time indices at
which agent $i$ computes its state updates\footnote{If computing a state update takes some non-zero number of timesteps, we can make $K^i$ the set of
times at which agent $i$'s computation of a state update completes. For simplicity we assume that computing a state update takes agent $i$ zero time and that state
updates are computed by agent $i$ at the points in time indexed by $K^i$.} by $K^i$, i.e., 
\begin{equation}
K \supseteq K^i := \{k \mid x^i_i \textnormal{ is updated by agent $i$ at time $k$}\}. 
\end{equation}
At times $k \in K \backslash K^i$ agent $i$ does not compute
any state updates and hence $x^i_i(k; t)$ does not change at these times, though $x^i_{-i}(k; t)$ 
can still change if a transmission
from another agent arrives at agent $i$ at time $k$. 
We note that $K$ and the sets $K^i$ need not be known by the agents as they are merely
tools used in the analysis of the forthcoming asynchronous algorithm. 
We also take $T = \N$ as the set of ticks of the dual update clock in the cloud
without loss of generality, though
there need not be any relationship between $T$ and $K$. 

Suppose that agent $j$ computes a state update at time $k_a$ and then begins transmitting
its state to agent $i$ also at time $k_a$. Due to communication delays, 
this transmission may not arrive at agent $i$ until, say, time $k_b > k_a$.
Suppose further that agent $j$'s next transmission to agent $i$ does not arrive
at agent $i$ until time $k_c > k_b$. It will
be useful in the following discussion to relate the time $k_b$ (at which the first
transmission arrives) 
to the time $k_a$ (at which it was originally computed by agent $j$).
Suppose at time $k$, with $\mu(t)$ onboard all agents, that agent $i$ has some value of agent $j$'s state,
denoted $x^i_j(k; t)$. 
We use $\tau^i_j(k)$ to denote
the time at which the value of $x^i_j(k; t)$ was originally computed by agent $j$. Above,
$\tau^i_j(k_b) = k_a$, and because the value of $x^i_j(\cdot; t)$ will not change again after $k_b$ until
time $k_c$, we have $\tau^i_j(k') = k_b$ for all $k_b \leq k' < k_c$. 
We similarly define $\tau^c_i : T \to K$ for all $i \in [N]$ to fulfill the same role for transmissions of state values
from agent $i$ to the cloud: at time $t$ in the cloud, $\tau^c_i(t)$ is the time $k$ at which
agent $i$ computed the state value it most recently sent to the cloud\footnote{The agents can send
multiple state values to the cloud between dual updates, though only the most recent transmission from agent $i$ to the cloud
will be kept by the cloud.}.
For all $i$, $j$, and $k$ we have $0 \leq \tau^i_j(k) \leq k$ by definition, and 
we impose the following assumption on $K^i$, $T$, $\tau^i_j$, and $\tau^c_i$
for all $i \in [N]$ and $j \in [N]$. 
\begin{assumption} \label{as:ttau}
For all $i \in [N]$ the set $K^i$ is infinite, and for a sequence $\{k_d\}_{d=1}^{\infty}$ in $K^i$ tending to infinity
we have
\begin{equation} \label{eq:atau}
\lim_{d \to \infty} \tau^i_j(k_d) = \infty
\end{equation}
for all $j \in \Ni_i$. Furthermore, the set $T$ is infinite and for a sequence $\{t_d\}_{d=1}^{\infty}$ in $T$
tending to infinity we have
\begin{equation} \label{eq:ctau}
\lim_{d \to \infty} \tau^c_i(t_d) = \infty
\end{equation}
for all $i \in [N]$. \hfill $\lozenge$
\end{assumption}

Requiring that $K^i$ be infinite guarantees that no agent will stop updating its state and Equation \eqref{eq:atau} guarantees that
no agent will stop sending state updates to its essential neighbors. 
Similarly, $T$ being infinite guarantees that the cloud continues to update $\mu$ and Equation \eqref{eq:ctau}
ensures that no agent stops sending its state to the cloud. 
Assumption~\ref{as:ttau} can therefore be understood
as ensuring that the system ``keeps running.'' 

For a fixed $\mu(t)$, agent $i$'s update law is written as follows (where $j \neq i$):
\begin{align}
x^i_i(k+1; t) &= \begin{cases} \Pi_{X_i}\Big[x^i_i(k; t) - \gamma\dxi \lve\big(x^i(k; t), \mu(t)\big)\Big] & k \in K^i \\
                                                                                               x^i_i(k; t) & k \not\in K^i \end{cases}\label{eq:agentiupdate} \\            
x^i_j(k+1; t) &= \begin{cases} x^j_j\big(\tau^i_j(k\!+\!1); t\big)                                         & \textnormal{$i$ receives $j$'s state at $k\!+\!1$} \\
                                                                                               x^i_j(k; t) & \textnormal{otherwise} \end{cases}\label{eq:agentij}. 
\end{align}
This update law has each agent performing gradient descent in its own state and waiting for other agents to update their states
and send them to the others in the network. 
This captures in a precise way that agent $i$ immediately incorporates transmissions from other agents into its current state value and that
such state values will generally be ``out-dated,'' as indicated by the presence of the $\tau^i_j$ term on the right-hand side of Equation \eqref{eq:agentij}.

\subsection{Cloud Dual Update Law}
While optimizing with $\mu(t)$ onboard, the agents compute some number of state updates using Equation~\eqref{eq:agentiupdate} and then send their states
to the cloud. It is not assumed that the agents send their states to the cloud at the same time or that they do so after the same number of state updates. 
Once the cloud has received states from the agents, it computes $\mu(t+1)$ and sends $\mu(t+1)$ to the agents, and then this process repeats. 
Assumption~\ref{as:ttau} specified that these operations do not cease being executed, and we impose the following basic assumption on the
sequence of updates that take place in the system. 
\begin{assumption} \label{as:agentcomms} 
\leavevmode
\begin{enumerate}[a.]
\item When the cloud sends $\mu(t+1)$ to the agents, it arrives in finite time. \label{as:ac2}
\item Any transmission originally sent from agent $i$ to agent $j$ while they have $\mu(t)$ onboard is only used by agent $j$ if it is received before $\mu(t+1)$.\label{as:ac3}
\item All transmissions arrive in the order in which they were sent\label{as:ac4}. 
\item \new{There is an increasing sequence of times $\{k_t\}_{t \in T}$ such that only $\mu(t)$ is used in the agents' state updates at timesteps $k \in K$
satisfying $k_t \leq k < k_{t+1}$.}
\label{as:ac5}
\hfill $\lozenge$
\end{enumerate}
\end{assumption}

Assumption~\ref{as:agentcomms}.\ref{as:ac2} is enforced simply to ensure
that the optimization process does not stall and is easily satisfied in practice. 
Assumption~\ref{as:agentcomms}.\ref{as:ac3} is enforced because $\mu(t)$ parameterizes
\begin{equation}
\hat{x}^t := \argmin_{x \in X} \lve(x, \mu(t)),
\end{equation}
which is the point the agents approach while optimizing with $\mu(t)$ onboard. 
Suppose that a message from agent $i$ is sent to agent $j$ while they have $\mu(t)$ onboard but is received after they have
$\mu(t+1)$ onboard. We will in general have $\hat{x}^t \neq \hat{x}^{t+1}$ so that the arrival of agent $i$'s message to agent $j$ effectively redirects
agent $j$ away from $\hat{x}^{t+1}$ and toward $\hat{x}^t$, delaying (or preventing) progress of the optimization algorithm. With respect to implementation,
little needs to be done to enforce this assumption. 
All communications between agents can be transmitted along with the timestamp $t$ of the dual
variable onboard the agent sending the message at the time it was sent. The agent receiving this message can compare the value of $t$ in the message
with the timestamp of the dual variable it currently has onboard, and any message with a mismatched timestamp can be discarded. 
Assumption~\ref{as:agentcomms}.\ref{as:ac3} can therefore be implemented in software without further constraining the 
agents' behavior or the optimization problem itself.

Assumption~\ref{as:agentcomms}.\ref{as:ac4} will be satisfied by a number of communication protocols, including
TCP \cite[Section 13]{comer00}, which is used on much of the internet, and does not constrain the agents' behavior because
it can be enforced in software by choosing an applicable communication protocol. 
It is enforced here
to prevent pathological behavior that can prevent the optimization from converging at all. 
\new{Assumption~\ref{as:agentcomms}.\ref{as:ac5} enforces that the agents use the same value of the dual
variable in their updates.} 
Assumption~\ref{as:agentcomms}.\ref{as:ac5} is the lone point of synchrony in the system
and is necessary for convergence of the asynchronous algorithm. This necessity is verified by a counter-example in Section~\ref{sec:counter} 
wherein violating only Assumption~\ref{as:agentcomms}.\ref{as:ac5} causes the system not to converge.

After the agents have taken some number of steps using $\mu(t)$ and have sent their states to the cloud, the cloud aggregates these states
into a vector which we denote $x^c_t$, defined as
\begin{equation}
x^c_t = \big(x^1_1\big(\tau^c_1(t); t\big), \ldots, x^N_N\big(\tau^c_N(t); t\big)\big).
\end{equation}
Then we adapt the dual update in Algorithm~\ref{alg:ensemble} to account for the time
the cloud spends waiting to receive transmissions from the agents, giving
\begin{equation} \label{eq:cloudupdate}
\mu(t+1) = \Pi_{\dualset}\left[\mu(t) + \rho\left(g(x^c_t) - \dreg\mu(t)\right)\right].
\end{equation}

\subsection{Asynchronous Primal-Dual Update Law} \label{ss:algo}
We now state the full asynchronous primal-dual algorithm that will be the focus of the remainder of the paper. 
Below we use the notation $C^i \subset K$ to denote the set of times at which agent $i$ sends its 
state to the cloud. \new{We also use the notation $\Rij$ to denote the set of times at which agent $j$
sends its state to agent $i$; if $j \not\in \Ni_i$, then $\Rij = \0$. 
Note that $C^i$ need not have any relationship to $K^i$ or $\Rij$ and that agent $i$ need not
know $C^i$ as it is merely a tool used for analysis. Similarly, $\Rij$ does not need
to have any relationship to $K^i$ or $C^i$ and does not need to be known by any agent.
We state the algorithm with the cloud waiting for each agent's state before computing a dual update
because this will typically be the desired behavior in a system. However, we do point out how to 
eliminate this assumption and the impact of this removal in Remark~\ref{rem:notallupdate} in the next section. }
\new{
\begin{talgorithm} \label{alg:async} 
$ $ \newline
Step 0: Initialize all agents and the cloud with $x(0) \in X$ and $\mu(0) \in \dualset$. Set $t = 0$ and $k = 0$. \\
Step 1: For all $i \in [N]$ and all $j \in \Ni_i$, 
if $k \in \Rij$, then agent $j$ sends $x^j_j(k; t)$ to agent $i$ (though it may not be received for some time). \\
Step 2: For all $i \in [N]$ and all $j \in \Ni_i$, execute
\begin{align}
x^i_i(k+1; t) &= \begin{cases} \Pi_{X_i}\left[x^i_i(k; t) - \gamma\dxi\lve\big(x^i(k; t), \mu(t)\big)\right] & k \in K^i \\
                                                                                                 x^i_i(k; t) & k \not\in K^i \end{cases} \\
x^i_j(k+1; t) &= \begin{cases} x^j_j\big(\tau^i_j(k+1); t\big)                                               & \textnormal{$i$ receives $j$'s state at time $k+1$} \\
                                                                                                 x^i_j(k; t) & \textnormal{otherwise} \end{cases}. 
\end{align}
Step 3: If $k+1 \in C^i$, agent $i$ sends $x^i_i(k+1; t)$ to the cloud. Set $k := k + 1$. 
If all components of $x^c_t$ have been updated since the agents received $\mu(t)$, the cloud computes
\begin{equation} 
\mu(t+1) = \Pi_{\dualset}\left[\mu(t) + \rho\left(g(x^c_t) - \dreg\mu(t)\right)\right]
\end{equation}
and sends $\mu(t+1)$ to the agents. Set $t := t + 1$. \\
Step 4: Return to Step 1. \hfill $\triangle$
\end{talgorithm}
} 

We show in Section~\ref{sec:main} that Algorithm~\ref{alg:async} approximately converges to $(\hxve, \hmve)$.

%% file: 5-convergence.tex
\section{Convergence of Asynchronous Primal-Dual Method} \label{sec:main}
In this section we examine the convergence properties of Algorithm~\ref{alg:async}
\new{and develop its convergence rates.} 
\new{
For clarity of presentation, we first show results that assume that
all agents send their states to the cloud before the cloud computes
each dual update. Once the main results of this section are established,
we explain how to eliminate this Assumption in Remark \ref{rem:notallupdate}. 
}

\subsection{Block Maximum Norm Basics}
First we consider the agents optimizing in the primal space with a fixed $\mu(t)$. 
We will examine convergence using a block-maximum norm similar to
that defined in Section 3.1.2 of \cite{bertsekas97}. First, for a vector $x \in X := X_1 \times \cdots \times X_N$,
we can decompose $x$ into its components as $x := (x_1, \ldots, x_N)$, and we refer to each
such component of $x$ as a \emph{block} of $x$; in Algorithm~\ref{alg:async}, agent $i$
updates block $i$ of $x^i$. Using the notion of a block we have the following
definition. 

\begin{definition} \label{def:bmax}
For a vector $x \in \R^n$ comprised of $N$ blocks, with the $i^{th}$
block being $x_i \in \R^{n_i}$,
the norm $\bmax{\cdot}$ is defined as\footnote{For concreteness we focus on the $\btup$-norm, though the results
we present can be extended to some more general weighted block-maximum norms of the form 
$\|x\|_{max} = \max_{i \in [N]} \|x_i\|_{p_i}/w_i$,
where ${w_i > 0}$ and ${p_i \in \N}$ for all $i \in [N]$.}
$\bmax{x} = \max_{i \in [N]} \|x_i\|_2$. \hfill $\Diamond$ 
\end{definition}

We have the following lemma regarding the matrix norm induced on $\R^{n \times n}$ by $\bmax{\cdot}$. 
In it, we use the notion of a block of a matrix. For $n = \sum_{i=1}^{N} n_i$,
the $i^{th}$ block of $A \in \R^{n \times n}$ is the $n_i \times n$ matrix formed by rows with indices
$\sum_{k=1}^{i-1} n_k + 1$ through $\sum_{k=1}^{i}n_k$ in $A$, and
we denote the $i^{th}$ block of $A$ by $A^{[i]}$.

\begin{lemma} \label{lem:matrixnorm}
For all $A \in \R^{n \times n}$,
\begin{equation}
\bmax{A} \leq \|A\|_2. 
\end{equation}
\end{lemma}
\emph{Proof:} For $A^{[i]}$ the $i^{th}$ block of $A$, let $A^{[i]}_{\ell,j}$ be the $\ell^{th}j^{th}$ entry
of that block and let $\Sp^{n-1}$ be the unit sphere in $\R^n$. Then for any $x \in \Sp^{n-1}$ and $i \in [N]$ we have
\begin{equation} \label{eq:block2norm}
\|A^{[i]}x\|_2 = \left(\sum_{k=1}^{n_i}\left(\sum_{j=1}^{n}A^{[i]}_{k,j}x_j\right)^2\right)^{1/2}.
\end{equation}
Each term on the right-hand side of Equation \eqref{eq:block2norm} is manifestly positive
so that summing over every block gives
\begin{equation}
\|A^{[i]}x\|_2 \leq \left(\sum_{i=1}^{N} \sum_{k=1}^{n_i} \left(\sum_{j=1}^{n} A^{[i]}_{k,j} x_j \right)^2 \right)^{1/2},
\end{equation}
which we can write in terms of $A$ as
\begin{equation}
\left(\sum_{i=1}^{N} \sum_{k=1}^{n_i} \left(\sum_{j=1}^{n} A^{[i]}_{k,j} x_j \right)^2 \right)^{1/2} = \left(\sum_{\ell=1}^{n} \left(\sum_{j=1}^{n} A_{\ell,j} x_j \right)^2 \right)^{1/2} = \|Ax\|_2.
\end{equation}
Then, taking the maximum over all blocks, we have
\begin{equation}
\bmax{Ax} = \max_{i \in [N]} \|A^{[i]}x\|_2 \leq \|Ax\|_2,
\end{equation}
for all $x \in \Sp^{n-1}$, and the result follows by taking
the supremum over $x \in \Sp^{n-1}$. 
\hfill $\blacksquare$

We also have the following elementary lemma relating norms of vectors.
\begin{lemma} \label{lem:vecnorm}
For all $x \in X$,
\begin{equation}
\|x\|_2^2 \leq N\bmax{x}^2 \quad \textnormal{ and } \quad \|x\|_2 \leq \sqrt{N}\bmax{x}.
\end{equation}
\end{lemma}
\emph{Proof:} 
We find
\begin{equation}
\|x\|_2^2 = \sum_{i \in [N]} \|x_i\|_2^2 \leq N\max_{i \in [N]} \|x_i\|_2^2 = N\bmax{x}^2,
\end{equation}
and then take the square root. \hfill $\blacksquare$

\subsection{Convergence in the Primal Space} \label{ss:sets}
We now examine what happens when the agents are optimizing in between transmissions from the cloud. 
We consider the case of some $\mu(t) \in \dualset$ fixed onboard all agents
and as in Section~\ref{sec:decent} we use
$\hat{x}^t = \argmin_{x \in X} \lve\big(x, \mu(t)\big)$. 
Under these conditions, the agents are
asynchronously minimizing the function
\begin{equation}
\pfunc(x) := \lve\big(x, \mu(t)\big)
\end{equation}
until $\mu(t+1)$ arrives from the cloud. 
To assess the convergence of Algorithm~\ref{alg:async} in the primal space, we define a sequence of
sets $\{X^t(s)\}_{s \in \N}$ to make use of the framework in Sections 6.1 and 6.2 of \cite{bertsekas97}. These sets must satisfy the following assumption
which is based on the assumptions in those sections. 

\begin{assumption} \label{as:xk}
With a fixed value of $t$ and a fixed $\mu(t)$ onboard all agents, the sets $\{X^t(s)\}_{s \in \N}$ satisfy:
\begin{enumerate}[a.]
\item $\cdots \subset X^t(s+1) \subset X^t(s) \subset \cdots \subset X$\label{as:xk1}
\item $\lim_{s \to \infty} X^t(s) = \{\hat{x}^t\}$\label{as:xk2} 
\item For all $i$, there are sets $X^t_{i}(s) \subset X_i$ satisfying\label{as:xk3}
\begin{equation}
X^t(s) = X^t_{1}(s) \times \cdots \times X^t_{N}(s) 
\end{equation}
\item For all $y \in X^t(s)$ and $i \in [N]$, $\updatei(y) \in X^t_{i}(s+1)$, where $\updatei(y) := \Pi_{X_i}\left[y_i - \gamma\dxi\pfunc(y)\right]$\label{as:xk4}.\hfill$\lozenge$
\end{enumerate} 
\end{assumption}

Unless otherwise noted, when writing $x \in X^t(s)$ for some vector $x$, the set $X^t(s)$ is chosen
with the largest value of $s$ that makes the statement true.
Assumptions~\ref{as:xk}.\ref{as:xk1} and \ref{as:xk}.\ref{as:xk2} require that we have a nested chain of sets to descend
that ends with $\hat{x}^t$. Assumption~\ref{as:xk}.\ref{as:xk3} allows the blocks of the primal
variable to be updated independently by the agents while still guaranteeing that progress
toward $\hat{x}^t$ is being made. Assumption~\ref{as:xk}.\ref{as:xk4} guarantees forward progress down the chain of
sets $\{X^t(s)\}_{s \in \N}$ whenever an agent computes a state update. More will be said about this assumption and its consequences
below in Remark~\ref{rem:no2step}. 

Recalling that $\plip$ is the (maximum, over $\mu \in M$) Lipschitz constant of $\dx \lve(\cdot, \mu)$, we define the constant
\begin{equation}
\qp = \max\{|1 - \gamma\preg|, |1 - \gamma\plip|\}. 
\end{equation}
We then have the following lemma that lets us determine the value of $\qp$ based upon $\gamma$.
\begin{lemma} \label{lem:q}
For $\gamma \in (0, 2/\plip)$ \new{and $\alpha \in (0, L_p)$} we have $\qp \in (0, 1)$. Furthermore the minimum value of $\qp$ is
\begin{equation}
\qp^* = \frac{\plip - \preg}{\plip + \preg} \quad \textnormal{ when } \quad \gamma = \frac{2}{\plip + \preg}.
\end{equation}
\end{lemma}
\emph{Proof:} See Theorem 3 on page 25 of \cite{poljak87}. \hfill $\blacksquare$

We proceed under the restrictions that $\gamma \in (0, 2/\plip)$ \new{and $\alpha \in (0, L_p)$ for the remainder of the paper}. 
\new{
To simplify the presentation of results in this section,
for all $t \in T$ we assume that all agents simultaneously received $\mu(t)$ at some time $k_t$. This $k_t$
serves the same role as in Assumption~\ref{as:agentcomms}.\ref{as:ac5}, though the agents do not actually need to
receive $\mu(t)$ at the exact same time and indeed any means of enforcing Assumption~\ref{as:agentcomms}.\ref{as:ac5}
will suffice. We retain this assumption on $k_t$ for the simplicity it provides below. 
}

For each $k_t$, we define the quantity
\begin{equation}
\dkt := \max_{i \in [N]} \bmax{x^i(k_t; t) - \hat{x}^t},
\end{equation}
which is the ``worst-performing'' block onboard any agent with respect to distance from $\hat{x}^t$. 
For a fixed value of $t$ we define each element in the sequence of sets $\{X^t(s)\}_{s \in \N}$ as
\begin{equation} \label{eq:xkdef}
X^t(s) = \left\{y \in X : \bmax{y - \hat{x}^t} \leq \qp^s\dkt\right\}.
\end{equation}
By definition, at time $k_t$ we have $x^i(k_t; t) \in X^t(0)$ for all $i$, and moving from $X^t(0)$ to
$X^t(1)$ requires contracting toward $\hat{x}^t$ (with respect to $\bmax{\cdot}$) by a factor of $\qp$. 
We have the following proposition.
\begin{proposition}
The collection of sets $\{X^t(s)\}_{s \in \N}$ as defined in Equation \eqref{eq:xkdef} satisfies Assumption~\ref{as:xk}. 
\end{proposition}
\emph{Proof:} By definition $X^t(s) \subset X$ for all $s$. From Equation \eqref{eq:xkdef}, we see that
\begin{equation}
X^t(s+1) = \left\{y \in X : \bmax{y - \hat{x}^t} \leq \qp^{s+1}\dkt\right\}.
\end{equation}
Because $\qp \in (0, 1)$, we have $\qp^{s+1} < \qp^s$, so that for $y \in X^t(s+1)$ we find
\begin{equation}
\bmax{y - \hat{x}^t} \leq \qp^{s+1}\dkt < \qp^s\dkt,
\end{equation} 
giving $y \in X^t(s)$ as well. Then $X^t(s+1) \subset X^t(s) \subset X$ for all $s \in \N$ and
Assumption~\ref{as:xk}.\ref{as:xk1} is satisfied. 

We see that
\begin{align}
\lim_{s \to \infty} X^t(s) &= \lim_{s \to \infty} \left\{y \in X : \bmax{y - \hat{x}^t} \leq \qp^s\dkt\right\} \\
                           &= \left\{y \in X : \bmax{y - \hat{x}^t} \leq 0\right\} \\ 
                           &= \{\hat{x}^t\},
\end{align}
which follows because $\bmax{\cdot}$ is a norm. Then Assumption~\ref{as:xk}.\ref{as:xk2} is satisfied as well. 

For Assumption~\ref{as:xk}.\ref{as:xk3}, the definition of $\bmax{\cdot}$ lets us easily decompose $X^t(s)$. In particular,
we see that 
\begin{equation}
\bmax{y - \hat{x}^t} \leq \qp^s\dkt \textnormal{ if and only if } \|y_i - \hat{x}^t_i\|_2 \leq \qp^s\dkt
\end{equation}
for all $i \in [N]$. Immediately then we have
\begin{equation}
X^t_{i}(s) = \left\{y_i \in X_i : \|y_i - \hat{x}^t_i\|_2 \leq \qp^s\dkt\right\}
\end{equation}
from which it is clear that
$
X^t(s) = X^t_{1}(s) \times \cdots \times X^t_{N}(s)
$
and thus that Assumption~\ref{as:xk}.\ref{as:xk3} is satisfied.

Finally, we show that Assumption~\ref{as:xk}.\ref{as:xk4} is satisfied. For a fixed $t$ and fixed $\mu(t)$ take some $y \in X^t(s)$. 
Recall the following exact expansion of $\dx\pfunc$:
\begin{align}
\dx\pfunc(y) - \dx\pfunc(\hat{x}^t) &= \int_{0}^{1} \dx^2\pfunc\big(\hat{x}^t + \tau(y - \hat{x}^t)\big) (y - \hat{x}^t) d\tau \\
                                    &= \left(\int_{0}^{1} \dx^2\pfunc\big(\hat{x}^t + \tau(y - \hat{x}^t)\big)d\tau\right) \cdot (y - \hat{x}^t) \\
                                    &=: \Hi(y)(y - \hat{x}^t)\label{eq:2taylor},
\end{align}
where we have defined $\Hi$ as 
\begin{equation} \label{eq:Hdef}
\Hi(y) := \int_{0}^{1} \dx^2\pfunc\big(\hat{x}^t + \tau(y - \hat{x}^t)\big)d\tau.
\end{equation}

Using the non-expansive property of $\Pi_{X_i}[\cdot]$ with respect to $\|\cdot\|_2$ on $\R^{n_i}$, we find
that for $y \in X^t(s)$
\begin{align}
\|\updatei(y) - \hat{x}_i^t\|_2 &=    \big\|\Pi_{X_i}\big[y_i - \gamma\dxi\pfunc(y)\big] - \Pi_{X_i}\big[\hat{x}_i^t - \gamma\dxi\pfunc(\hat{x}^t)\big]\big\|_2 \\
                                &\leq \|y_i - \gamma\dxi\pfunc(y) - \hat{x}_i^t + \gamma\dxi\pfunc(\hat{x}^t)\|_2 \\
                                &\leq \max_{i \in [N]} \|y_i - \gamma\dxi\pfunc(y) - \hat{x}_i^t + \gamma\dxi\pfunc(\hat{x}^t)\|_2 \\
                                &= \bmax{y - \hat{x}^t - \gamma\left(\dx \pfunc(y) - \dx\pfunc(\hat{x}^t)\right)} \\
                                &= \bmax{y - \hat{x}^t - \gamma\Hi(y)(y - \hat{x}^t)} \\
                                &\leq \bmax{I - \gamma\Hi(y)} \bmax{y - \hat{x}^t} \\
                                &\leq \|I - \gamma\Hi(y)\|_2 \bmax{y - \hat{x}^t},\label{eq:almostdone} 
\end{align} 
where the third equality follows from Equation \eqref{eq:2taylor} and where the last inequality follows from Lemma~\ref{lem:matrixnorm}. 
For any choice of $\mu \in M$, the $\preg$-strong monotonicity and $\plip$-Lipschitz properties of $\dx \lve(\cdot, \mu(t))$ give
$\preg I \preceq \Hi(\cdot) \preceq \plip I$,
which implies that the eigenvalues
of $\Hi(\cdot)$ are bounded above by $\plip$ and below by $\preg$
for all $\mu(t) \in M$. 
Using this fact and that $\Hi(y)$ is symmetric, we see that
\begin{align}
\|I - \gamma\Hi(y)\|_2 &= \max\Big\{|\lambda_{min}(I - \gamma\Hi(y))|, |\lambda_{max}(I - \gamma\Hi(y)|\Big\} \\
                       &= \max\{|1 - \gamma\preg|, |1 - \gamma\plip|\} \\
                       &= \qp, 
\end{align}
where $\lambda_{min}$ and $\lambda_{max}$ denote the minimum and maximum eigenvalues of a matrix, respectively. 
Then from Equation \eqref{eq:almostdone}, and the fact that $\bmax{y - \hat{x}^t} \leq \qp^s\dkt$ by hypothesis, we find
\begin{equation}
\|\updatei(y) - \hat{x}_i^t\|_2 \leq \qp\bmax{y - \hat{x}^t} \leq \qp^{s+1}\dkt,
\end{equation}
so that $\updatei(y) \in X^t_i(s+1)$ as desired. \hfill $\blacksquare$

We comment on one consequence of Assumption~\ref{as:xk} in particular below where we use
the notation
\begin{equation}
X^t_{-i}(s) : = X^t_1(s) \times \cdots \times X^t_{i-1}(s) \times X^t_{i+1}(x) \times \cdots \times X^t_N(s). 
\end{equation}

\begin{remark} \label{rem:no2step}
Suppose at time $k$
agent $i$ has state vector $x^i(k; t) \in X^t(s)$ and $k+1 \in K^i$. Then Assumption~\ref{as:xk}.\ref{as:xk4} implies that
${x^i_i(k+1; t) = \updatei(x^i(k; t)) \in X^t_{i}(s+1)}$. Suppose, before any other agent transmits an updated state to agent $i$, that
agent $i$ performs another update of its own state. Just before the second update, 
$x^i(k+1; t)$ is equal to $x^i(k; t)$ with the entry for $x^i_i$ replaced with the update just computed, $\updatei(x^i(k; t))$;
all other entries of $x^i(k+1; t)$ remain unchanged from $x^i(k; t)$. Because no
other agents' states have changed, we still have $x^i_{-i}(k+1; t) \in X^t_{ -i}(s)$ and, as a result, $x^i(k+1; t) \in X^{t}(s)$.
In general, $x^i(k+1; t) \not\in X^t(s+1)$ here precisely because $x^i_{-i}(k+1; t)$ has not changed. Due to the fact that $x^i(k+1; t) \in X^t(s)$,
the second update performed by agent $i$ results in $\updatei(x^i(k+1; t)) \in X^t_{i}(s+1)$ once more, though, in general,
no further progress, e.g., to $X^t_i(s+2)$, can be made without further updates from the other agents. Then while an agent is waiting
for updates from other agents, its progress toward $\hat{x}^t$ can be halted, though it does not ``regress'' backwards from, say, $X^t_{i}(s)$ to
$X^t_{i}(s-1)$. 
\hfill $\Diamond$
\end{remark}

We proceed to use Assumption~\ref{as:xk} to estimate the primal convergence rate of Algorithm~\ref{alg:async}.

\subsection{Single-Cycle Primal Convergence Rate Estimate}
The structure of the sets $\{X^t(s)\}_{s \in \N}$ enables us to extract a convergence rate estimate in the primal space. To demonstrate this point, consider
the starting point of the algorithm: all agents have onboard some state $x(0) \in X^0(0)$ and dual vector $\mu(0) \in \dualset$. 
Suppose that agent $i$ takes a single gradient descent step, say at time $k^i \in K^i$, from $x(0)$ with $\mu(0)$ held fixed. 
From Assumption~\ref{as:xk}.\ref{as:xk4}, this results in
agent $i$ having $x^i_i(k^i; 0) = \updatei\big(x^i(0; 0)\big) = \updatei\big(x(0)\big)\in X^0_{i}(1)$. Once agent $i$ transmits the state $x^i_i(k^i; 0)$ to 
its essential neighbors, and once all other
agents themselves have taken a descent step and
transmitted their states to their essential neighbors, say, at time $\bar{k}$, agent $i$ will have $x^i(\bar{k}; 0) \in X^{0}(1)$.
Agent $i$'s next descent step, say at time $\ell > \bar{k}$, then results in $x^i_i(\ell; 0) = \updatei\big(x^i(\bar{k}; 0)\big) \in X^0_{i}(2)$. Then the process of communicating and descending
repeats. 
To keep track of how many times this process repeats, we have the following definition.

\begin{definition} \label{def:cycles}
After the agents all have just received $\mu(t)$, the first \emph{cycle} ends as soon as (i) each agent has computed a state update
and (ii) each agent has sent that updated state to all of its essential neighbors and it has been received by them. Subsequent cycles are completed
when the preceding cycle has ended and criteria (i) and (ii) are met again, with any number of cycles possible
between $k_t$ and $k_{t+1}$. \hfill $\Diamond$
\end{definition}

It is possible for one agent to compute and share several state updates with the other agents within
one cycle if some other agent in the network is updating
more slowly. 
For a fixed value of $\mu(0)$ onboard all agents and a common initial state $x(0)$, the first cycle 
will move each agent's copy of the ensemble state from $X^0(0)$ to $X^0(1)$, the second cycle
will move it from $X^0(1)$ to $X^0(2)$, etc. When the agents have $\mu(t)$ onboard,
we use $c(t)$ to denote the number of cycles the agents complete before the first agent sends its state to the cloud for use in computing $\mu(t+1)$. 
We see that with $\mu(0)$ onboard, the agents complete $c(0)$ cycles to reach the set $X^0(c(0))$ and the cloud
therefore uses an element of $X^0(c(0))$ to compute $\mu(1)$. 
In particular, using Assumption~\ref{as:xk}.\ref{as:xk4} and the construction of the sets $\{X^0(s)\}_{s \in \N}$, 
this means that the convergence rate is geometric in the number of cycles completed:
\begin{equation}
\bmax{x^c_0 - \hat{x}^0} \leq \qp^{c(0)}D(k_0) = \qp^{c(0)}\bmax{x(0) - \hat{x}^0}. 
\end{equation}

Crucially, it need not be the case that all agents have the same state at the beginning of each cycle for this rate estimate to apply.
We show this in deriving a general primal convergence rate in the following lemma. 

\begin{lemma} \label{lem:primalconv}
Let Assumptions~\ref{as:sets}-\ref{as:xk} hold, let $\reg = (\preg, \dreg)$ be fixed, and let $\gamma \in (0, 2/\plip)$. 
When the agents are all optimizing with $\mu(t)$ onboard and the first agent sends its state to the cloud after $c(t)$ cycles,
we have
\begin{equation} \label{eq:pconvlast}
\bmax{\cloudxt - \hat{x}^t} \leq q_p^{c(t)}\dkt. 
\end{equation}
\end{lemma}
\emph{Proof:} 
Suppose the agents just received $\mu(t)$ from the cloud. For all $i \in [N]$ we have $x^i(k_t; t) \in X^t(0)$
from the definition of $\dkt$.
Then when agent $i$ computes a state update the result is $\updatei\big(x^i(k_t; t)\big) \in X^t_i(1)$.
When the agents have completed one cycle after receiving $\mu(t)$,
say by time $\bar{k}$, we have $x^i(\bar{k}; t) \in X^t(1)$. 
Iterating this process,
after $c(t)$ cycles agent $i$'s copy of the ensemble state moves from $X^t(0)$ to $X^t(c(t))$ for all $i \in [N]$. Then, 
when the agents send their states to the cloud, agent $i$ sends an element of $X^t_i(c(t))$. 
Then we have $x^c_t \in X^t(c(t))$ and, by definition,
$\bmax{x^c_t - \hat{x}^t} \leq q_p^{c(t)}\dkt$. \hfill $\blacksquare$ 

One can impose further assumptions, e.g., that a cycle occurs every $B$ ticks of $K$, in which 
case the exponent of $q_p$ in Equation \eqref{eq:pconvlast} becomes $\floor{(k_{t+1} - k_t)/B}$,
though for generality we do not do so. 

\subsection{Overall Convergence Rates}
Towards providing a convergence rate estimate in the dual space, we first present the following lemma. 

\begin{lemma} \label{lem:coco}
For any primal-dual pairs $(x_1, \mu_1) \in X \times M$ and $(x_2, \mu_2) \in X \times M$ such that
\begin{equation}
x_1 = \argmin_{x \in X} \lve(x, \mu_1) \textnormal{ and } x_2 = \argmin_{x \in X} \lve(x, \mu_2)
\end{equation}
we have
\begin{equation}
(\mu_2 - \mu_1)^T(g(x_1) - g(x_2)) \geq \frac{\preg}{\gradgm^2} \|g(x_1) - g(x_2)\|^2. 
\end{equation}
In addition,
\begin{equation}
\|\mu_1 - \mu_2\| \geq \frac{\preg}{\gradgm}\|x_1 - x_2\|. 
\end{equation}
\end{lemma}
\emph{Proof:} See \cite[Lemma 4.1]{koshal11}. \hfill $\blacksquare$

We now prove approximate convergence in the dual space and estimate the rate of convergence
there. 

\begin{theorem} \label{thm:dualconv}
Let all hypotheses of Lemma~\ref{lem:primalconv} hold and let the dual step-size satisfy
\begin{equation}
0 < \rho < \rho_0 := \min\left\{\frac{2\preg}{\gradgm^2 + 2\preg\dreg}, \frac{2\dreg}{1 + \dreg^2}\right\},
\end{equation}
where $M_g = \max_{x \in X} \|\nabla g(x)\|$. 
Then for all $t \geq 0$ 
\begin{equation}
\|\mu(t+1) - \hmve\|^2 \leq q_d^{t+1}\|\mu(0) - \hmve\|^2 + \sum_{\ell=0}^{t} q_d^{t-\ell}\Big(q_dN\gradgm^2\xdistmax^2q_p^{2c(\ell)} + 2\sqrt{N}\rho^2\gradgm^2\xdistmax\xdiam q_p^{c(\ell)}\Big),
\end{equation}
where $(0, 1) \ni \qd := (1 - \rho\dreg)^2 + \rho^2$, $\xdiam := \max_{x, y \in X} \|x - y\|$ 
is the diameter of $X$, and $\xdistmax := \max_{i \in [N]} \max_{x_i, y_i \in X_i} \|x_i - y_i\|$ is the maximum diameter
among the sets $X_i$, $i \in [N]$. 
\end{theorem}
\emph{Proof:} 
Using the non-expansive property of the projection operator $\Pi_M[\cdot]$ and expanding we find
\begin{align}
\|\mu(t+1) - \hmve\|^2 &=  \big\|\Pi_M\left[\mu(t) + \rho\left(g(x^c_t) - \dreg\mu(t)\right)\right] - \Pi_M\left[\hmve + \rho\left(g(\hxve) - \dreg\hmve\right)\right]\big\|^2 \\
\begin{split}
            &\leq (1 - \rho\dreg)^2\|\mu(t) - \hmve\|^2 + \rho^2\|g(\hxve) - g(x^c_t)\|^2 \\
            &\qquad\qquad- 2\rho(1 - \rho\dreg)\big(\mu(t) - \hmve\big)^T\big(g(\hxve) - g(x^c_t)\big). 
\end{split}\\
\end{align}
Adding $g(\hat{x}^t) - g(\hat{x}^t)$ inside the last set of parentheses, expanding, and 
applying Lemma~\ref{lem:coco} then gives
\begin{multline} \label{eq:goback1}
\|\mu(t+1) - \hmve\|^2 \leq (1 - \rho\dreg)^2\|\mu(t) - \hmve\|^2 + \rho^2\|g(\hxve) - g(x^c_t)\|^2 \\ - 2\rho(1 - \rho\dreg)\frac{\preg}{\gradgm^2}\|g(\hxve) - g(\hat{x}^t)\|^2 
                            -2\rho(1 - \rho\dreg)\big(\mu(t) - \hmve\big)^T\big(g(\hat{x}^t) - g(x^c_t)\big).
\end{multline}
Next, we have
\begin{equation}
0 \leq \|(1 - \rho\dreg)(g(\hat{x}^t) - g(x^c_t)) + \rho(\mu(t) - \hmve)\|^2,
\end{equation}
where expanding and re-arranging gives
\begin{equation}\label{eq:normeq}
{-2}\rho(1 - \rho\dreg)(\mu(t) - \hmve)^T(g(\hat{x}^t) - g(x^c_t)) \leq (1 - \rho\dreg)^2\|g(\hat{x}^t) - g(x^c_t)\|^2 + \rho^2\|\mu(t) - \hmve\|^2.
\end{equation}
Substituting Equation \eqref{eq:normeq} into Equation \eqref{eq:goback1} then gives
\begin{multline} \label{eq:goback2}
\|\mu(t+1) - \hmve\|^2 \leq \big((1 - \rho\dreg)^2 + \rho^2\big)\|\mu(t) - \hmve\|^2 + \rho^2\|g(\hxve) - g(x^c_t)\|^2 \\
     {-2}\rho(1 - \rho\dreg)\frac{\preg}{\gradgm^2}\|g(\hxve) - g(\hat{x}^t)\|^2 
     +(1 - \rho\dreg)^2\|g(\hat{x}^t) - g(x^c_t)\|^2.
\end{multline}

Next, we see that
\begin{align} 
\|g(\hxve) - g(x^c_t)\|^2 &= \|g(\hxve) - g(\hat{x}^t) + g(\hat{x}^t) - g(x^c_t)\|^2 \\
                          &\leq \|g(\hxve) - g(\hat{x}^t)\|^2 + \|g(\hat{x}^t) - g(x^c_t)\|^2 + 2\|g(\hxve) - g(\hat{x}^t)\|\|g(\hat{x}^t) - g(x^c_t)\|\label{eq:normeq2},
\end{align}
and substituting Equation \eqref{eq:normeq2} into Equation \eqref{eq:goback2} gives
\begin{align} \label{eq:goback3}
\|\mu(t+1) - \hmve\|^2 \leq &\big((1 - \rho\dreg)^2 + \rho^2\big)\|\mu(t) - \hmve\|^2 +\left(\rho^2 -2\rho(1 - \rho\dreg)\frac{\preg}{\gradgm^2}\right)\|g(\hxve) - g(\hat{x}^t)\|^2 \\
     &\quad+\big((1 - \rho\dreg)^2 + \rho^2\big)\|g(\hat{x}^t) - g(x^c_t)\|^2 +2\rho^2\|g(\hxve) - g(\hat{x}^t)\|\|g(\hat{x}^t) - g(x^c_t)\|.
\end{align}

Taking $\rho \in (0, \rho_0)$, we find that
\begin{equation} \label{eq:r1neg}
\rho^2 - 2\rho(1 - \rho\dreg)\frac{\preg}{\gradgm^2} < 0 \textnormal{ and } q_d \in (0, 1). 
\end{equation}
Substituting Equation \eqref{eq:r1neg} into Equation \eqref{eq:goback3}, 
and using the Lipschitz property of $g$, we find
\begin{equation} \label{eq:t1nearend}
\|\mu(t+1) - \hmve\|^2 \leq q_d\|\mu(t) - \hmve\|^2 + q_d\gradgm^2\|\hat{x}^t - x^c_t\|^2 +2\rho^2\gradgm^2\xdiam\|\hat{x}^t - x^c_t\|.
\end{equation} 
Lemmas~\ref{lem:vecnorm} and \ref{lem:primalconv} imply that 
\begin{equation}
\|\hat{x}^t - x^c_t\|_2 \leq \sqrt{N}\bmax{\hat{x}^t - x^c_t} \leq \sqrt{N} q_p^{c(t)}\xdistmax.
\end{equation}
Using this in Equation \eqref{eq:t1nearend} gives
\begin{equation}
\|\mu(t+1) - \hmve\|^2 \leq q_d\|\mu(t) - \hmve\|^2 + q_dN\gradgm^2\xdistmax^2q_p^{2c(t)} + 2\sqrt{N}\rho^2\gradgm^2\xdistmax\xdiam q_p^{c(t)},
\end{equation} 
where the result follows by summing over $t$. 
\hfill $\blacksquare$

\begin{remark}
Theorem~\ref{thm:dualconv} shows that convergence in the dual space is governed by three terms. The
first term decays as $q_d^{t+1}$ and represents a contraction toward $\hmve$. The next two terms
are essentially error terms that result from $x^c_t$ not equaling $\hat{x}^t$; to see this, note
that an exact dual method would have $c(t) = \infty$, causing the sum in Theorem~\ref{thm:dualconv}
to vanish, leaving only the contracting term. We see that larger values
of $c(t)$ lead $x^c_t$ closer to $\hat{x}^t$, causing the algorithm to approximate an ordinary
dual algorithm, thus leading to smaller errors. 

In addition, the $q_d^{t-\ell}$ term outside the sum
indicates that past errors contribute less to the overall dual error, with old error terms
accumulating powers of $q_d$ over time. To make faster progress using the asynchronous algorithm, one can 
have the agents perform small numbers of cycles for small values of $t$ and then increase
$c(t)$ as $t$ becomes large. Such a strategy makes later error terms small while weighting
earlier error terms only minimally, giving a small overall error. \hfill $\Diamond$
\end{remark}

We now present a result on primal convergence in Algorithm~\ref{alg:async}. 

\begin{theorem} \label{thm:pconv}
Let the primal step-size $\gamma \in (0, 2/\plip)$ and let all hypotheses of Lemma~\ref{lem:primalconv}
hold. Then for the sequence of primal vectors aggregated by the
cloud, $\{x^c_t\}_{t \in \N}$, we have
\begin{equation}
\|x^c_t - \hxve\|_2 \leq q_p^{c(t)}\sqrt{N}\xdistmax + \frac{\gradgm}{\preg}\|\mu(t) - \hmve\|_2. 
\end{equation}
\end{theorem}
\emph{Proof}: 
Adding $\hat{x}^t - \hat{x}^t$ and using Lemmas~\ref{lem:vecnorm}, \ref{lem:primalconv}, and \ref{lem:coco} we find
\begin{align}
\|x^c_t - \hxve\| &\leq \|x^c_t - \hat{x}^t\| + \|\hat{x}^t - \hxve\| \\
                  &\leq \sqrt{N}\qp^{c(t)}\xdistmax + \frac{\gradgm}{\preg}\|\mu(t) - \hmve\|,
\end{align}
where we have bounded $D(k_t)$ by $\xdistmax$. 
\hfill $\blacksquare$

Convergence in the primal space is then governed by two terms, one of which behaves
like a contraction whose exponent is $c(t)$ and the other which is a constant multiple
of the dual error, and we again find that completing more cycles improves accuracy. 
\new{
We also have
the following tradeoff between speed and accuracy induced by the
regularization of $L$. 

\begin{remark} \label{rem:tradeoff}
Theorem~\ref{thm:pconv}, Proposition~\ref{prop:regchoose} and 
Lemmas~\ref{lem:q} and \ref{lem:primalconv} together reveal a fundamental tradeoff between
convergence rate and accuracy in the primal space. On the one hand, Proposition~\ref{prop:regchoose} shows
that smaller values of $\preg$ lead to smaller errors while larger values of
$\preg$ lead to larger errors. On the other hand, Lemma~\ref{lem:q}
shows that larger values of $\preg$ lead to smaller values of $\qp$, and both Lemma~\ref{lem:primalconv}
and Theorem~\ref{thm:pconv}
show that smaller values of $\qp$ lead to faster convergence through the primal space, while
smaller values of $\preg$ cause $\qp$ to approach the value $1$, thereby slowing convergence. 
Then smaller values of $\preg$ lead to smaller errors at the expense of slower convergence, while
larger values of $\preg$ cause the system to converge more quickly, but to a point that is further
away from $(\hxve, \hmve)$. 

A similar tradeoff applies to $\dreg$ as well: Lemmas~\ref{lem:costdiff} 
and \ref{lem:gviolation} show that smaller values of $\beta$ can lead to smaller
errors, though the definition of $\qd$ in Theorem~\ref{thm:dualconv} shows that
a larger value of $\dreg$ decreases $\qd$, leading to faster convergence. The appropriate
balance of convergence speed and accuracy of a solution depends upon the problem
being solved, though, taken together, these results give one the tools to quantitatively
balance these two objectives. 

One can also see the use of regularizing $L$ in Theorems~\ref{thm:dualconv} and \ref{thm:pconv}.
If one were to set $\preg = 0$, then we would find $\qp = 1$ and primal updates would not
make any progress toward $\hat{x}^t$ in Lemma~\ref{lem:primalconv}. Such a case would
also cause the construction of the sets $\{X^t(s)\}_{s \in \N}$ to break down as no
``descent'' down this sequence of sets could be shown. Similarly, if one were to set
$\dreg = 0$, we would find $\qd = 1 + \rho^2$, in which case the only way to avoid
moving \emph{away} from $\hmve$ in the dual space would be to set $\rho = 0$,
thereby forestalling all progress in the dual space. Through their roles
in determining $\qp$ and $\qd$ (and the use of these constants 
 in the convergence analysis presented), it is evident that regularizing
with $\preg$ and $\dreg$ is essential to the analysis presented here.
\hfill $\Diamond$
\end{remark}
}

\new{
We now point out how to formulate convergence rate estimates without having
each agent send a state update to the cloud before it computes each dual
update.
\begin{remark} \label{rem:notallupdate}
If one allows the cloud to compute dual updates before receiving a state update
from each agent, then Lemma~\ref{lem:primalconv} should be modified to account
for only some values of $x^c_t$ changing from $x^c_{t-1}$. In particular, if
$N(t)$ agents send state updates to the cloud before it computes
$\mu(t+1)$ and $M(t) := N - N(t)$ do not, we find
\begin{equation}
\|x^c_t - \hat{x}^t\|_2 \leq \sqrt{N(t)\qp^{c(t)}\dkt + M(t)\xdistmax}.
\end{equation}
Propagating this through Theorems~\ref{thm:dualconv} and \ref{thm:pconv} gives
overall primal and dual convergence rate estimates for this case as well. 
In doing so, one finds that
executing a cloud update without a state update from each agent
can significantly harm convergence and it will usually be preferred to have the cloud
wait until it has received state information from all agents before each dual update. 
\hfill $\Diamond$
\end{remark}
}

%% file: 6-counterexample.tex
\section{Non-Convergence of the Asynchronous Dual Case} \label{sec:counter}

In this section we provide a counterexample to show that Assumption~\ref{as:agentcomms}.\ref{as:ac5} is 
necessary for the convergence of Algorithm~\ref{alg:async}. In it, we allow the agents to have
different values of the system's dual variable
and show that these differences can cause the primal and dual trajectories in Algorithm~\ref{alg:async} 
not to converge at all.
As will be shown, this is true even when each agent receives the most recent dual value 
at regular intervals and when the agents keep their states synchronized at all times. 

\setcounter{algorithm}{2}
\begin{algorithm}[!tp]
\caption{Asynchronous Dual Counterexample} \label{alg:counter}
\begin{algorithmic}[1]
\State Initialize $\mu \gets 0$, $\mu_{old} \gets 0$, $x_1 \gets 0$, and $x_2 \gets 0$.
\For{$\tau_{outer} = 1 \textrm{ to } 10$}
\For{$\tau_1 = 1 \textrm{ to } 500$} \% \emph{Mode 1}
    \State $x_2 \gets \theta_2(x_1, x_2, \mu_{old})$
    \State $x_1 \gets \theta_1(x_1, x_2, \mu)$
    \State $\mu \gets \theta_M(x_1, x_2, \mu)$    
\EndFor
\State $\mu_{old} \gets \mu$ 
\For{$\tau_2 = 1 \textrm{ to } 1500$} \% \emph{Mode 2}
    \While{$|x_1 - \theta_1(x_1, x_2, \mu)| > 10^{-5}$} 
    \State $x_1 \gets \theta_1(x_1, x_2, \mu)$
    \EndWhile
    \While{$|x_2 - \theta_2(x_1, x_2, \mu_{old})| > 10^{-5}$}
    \State $x_2 \gets \theta_2(x_1, x_2, \mu_{old})$
    \EndWhile
    \State $\mu \gets \theta_M(x_1, x_2, \mu)$
\EndFor
\State $\mu_{old} \gets \mu$
\EndFor
\end{algorithmic}
\end{algorithm}

The problem consists of two agents with scalar states and per-agent objectives
$f_1(x_1) = 0.1x_1$ and $f_2(x_2) = -0.1x_2$,
coupling cost $c \equiv 0$, and the constraint
${g(x) = \frac{1}{2}(x_1 - x_2)^2 - 0.2 \leq 0}$.
The regularization parameters were chosen to be $\preg = \dreg = 0.01$ and
the constraint set is $X = [0, 5]^2$. 

\begin{figure}[!tp]
  \centering
  \includegraphics[width=3in]{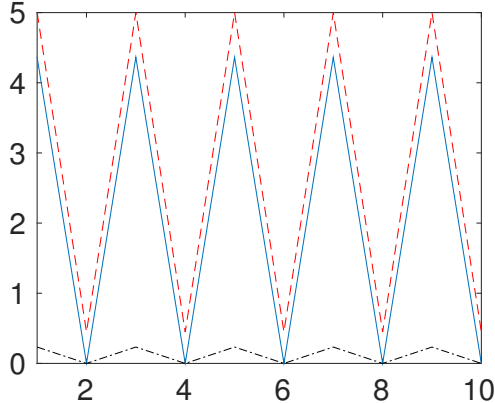}
  \caption{Primal and dual trajectories resulting from a simulation of Algorithm~\ref{alg:counter}, with $x_1$ the upper solid line,
  $x_2$ the upper dashed line, and $\mu$ the lower dash-dotted line.
  These oscillations are of constant magnitude and do not decay. 
  All terms are plotted at the end of each iteration of the outer loop.}
  \label{fig:counterplot}
\end{figure}

In this example, we will sometimes have one agent using an old dual value for some period
of time, and we denote this value by $\mu_{old}$; its value only changes when we write $\mu_{old} \gets \mu$
in the pseudocode in Algorithm~\ref{alg:counter}.
Otherwise, $\mu_{old}$ does not update with $\mu$. Similarly, the values of $x_1$ and $x_2$ only change when
explicitly updated below and operations listed sequentially below
actually occur sequentially so that the agents are updating at different times. 
To highlight the impact of asynchrony in the dual variable,
each agent always uses the most recent state of the other agent in its computations,
i.e., $x^1_2 = x^2_2$ and $x^2_1 = x^1_1$. Then there is no disagreement about state
values in the network and superscript indices are therefore omitted. 
For clarity, we write each argument of $\theta_1$ and $\theta_2$ out explicitly, including
specifying which dual variable is being used. To simplify notation, timestamps
are omitted in the pseudocode in Algorithm~\ref{alg:counter}. 

We have $L_p \approx 50.014$ so that using the stepsize bounds in Theorems~\ref{thm:dualconv} and
\ref{thm:pconv} we select $\gamma = 0.002$ and $\rho = 0.0003$. 
This example consists of alternating between two modes, shown in Algorithm~\ref{alg:counter} where
the dual update law in the cloud is represented by the symbol $\theta_M$.

Oscillations are shown in Figure~\ref{fig:counterplot} where we plot the primal and dual trajectories of a simulation 
implementing Algorithm~\ref{alg:counter}. Both states and the dual variable oscillate in a 
non-decaying fashion, indicating that Algorithm~\ref{alg:async} is not converging at all. 
We note here that synchronizing the dual variable in this example does indeed lead to convergence, indicating
that the asynchrony of the dual values is the source of oscillations and that Assumption~\ref{as:agentcomms}.\ref{as:ac5}
is a necessary condition for convergence in Algorithm~\ref{alg:async}.

%% file: 7-simulation.tex
\section{Simulation Results} \label{sec:simulation}
\begin{figure}[!tp]
\centering
\includegraphics[width=3in]{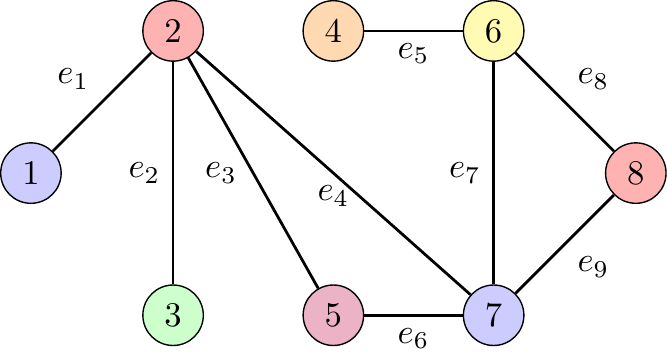}
\caption{The network across which $8$ agents route
traffic. There are $9$ edges, each with a maximum
capacity of $10$, and $8$ nodes. The edges used by each agent are listed in Table~\ref{tab:paths}.}
\label{fig:network}
\end{figure}

\begin{table}[!tp]
\centering
\begin{tabular}{|c|c|c|} \hline
Agent Number & Start Node$\to$End Node & Edges Traversed \\ \hline
$1$ & $1 \to 7$ & $e_1$, $e_3$, $e_6$ \\ \hline
$2$ & $2 \to 8$ & $e_4$, $e_7$, $e_8$ \\ \hline
$3$ & $3 \to 4$ & $e_2$, $e_4$, $e_7$, $e_5$ \\ \hline
$4$ & $5 \to 6$ & $e_3$, $e_4$, $e_7$ \\ \hline
$5$ & $1 \to 4$ & $e_1$, $e_3$, $e_6$, $e_7$, $e_5$ \\ \hline
$6$ & $3 \to 8$ & $e_2$, $e_4$, $e_9$ \\ \hline
$7$ & $4 \to 5$ & $e_5$, $e_8$, $e_9$, $e_6$ \\ \hline
$8$ & $6 \to 2$ & $e_7$, $e_4$ \\ \hline
\end{tabular} 
\caption{The edges traversed by each agent's flow.}
\label{tab:paths}
\end{table}
We now present simulation results for Algorithm~\ref{alg:async}.
We first discuss the problem to be solved and then cover 
our implementation. We then present numerical results that demonstrate convergence
of Algorithm~\ref{alg:async} on the cloud-based system and the
tradeoff between convergence rate and accuracy that is induced
by the Tikhonov regularization of $L$. 

\subsection{Problem Overview}
We consider a problem of routing $N = 8$ flows through a network
consisting of $8$ nodes and $9$ edges, representing, e.g., traffic
flow or sending data across a communication network, 
and each agent's decision variable is the flow rate of its
data through the network,
which is depicted in Figure~\ref{fig:network}. The nodes of the network
are not the agents themselves, but, instead, the agents are users
of the network attempting to route
traffic between certain pairs of these nodes.
The starting points, ending points, and edges which comprise the path traversed
by each flow are listed in Table~\ref{tab:paths}. 

We define the set
$\mathcal{E} := [9]$
to be the indices of the edges in the network. 
The cost of each agent is
$f_i(x_i) = -\delta_i\log(1 + x_i)$,
and we have selected $\delta_i = 100$ for all $i \in [8]$. 
The network also has an associated congestion cost
$c(x) = \frac{1}{20}x^TA^TAx$,
where 
\begin{equation}
A_{k,i} = \begin{cases} 1 & \textnormal{ if flow $i$ traverses edge $k$ } \\
                        0 & \textnormal{ otherwise}\end{cases}
\end{equation}
defines the network's adjacency matrix $A$.

Each edge in the network is subject to capacity constraints, expressed by requiring
$Ax \leq b$, where $b_i = 10$ for all $i \in \mathcal{E}$.
In addition, each flow rate is confined to $[0, 10]$, giving
$X = [0, 10]^8$. 
To demonstrate the effects of different values of $\preg$ and $\dreg$, 
three simulation runs were run: the first with $\preg = \dreg = 0.1$,
the second with $\preg = \dreg = 0.01$, and the third with $\preg = \dreg = 0.001$.
By sweeping $\preg$ and $\dreg$
across three orders of magnitude, we demonstrate the speed-accuracy tradeoff
discussed in Remark~\ref{rem:tradeoff}. 
For each $\preg$, we take
$\gamma = 2/(L_p + \preg)$,
and we take $\rho = 0.9\rho_0$ for each $(\preg, \dreg)$ pair. 

\subsection{Implementation and Numerical Results}
The implementation of the above problem allowed
as many quantities as possible to be random to demonstrate
asynchronous behavior. The
time between cloud updates was a random integer chosen
from the range $5$ to $100$ (inclusive) with uniform probability, and this number
represents the number of ticks of the virtual clock $K$ between $k_t$
and $k_{t+1}$. At each tick of $K$, each agent computed a state
update with probability $p_{update} = 0.05$ for all agents. 

The
communication graph at each tick of $K$ was an Erd\H{o}s-R\'{e}nyi graph \cite[Chapter 5]{mesbahi10},
which is a random graph wherein each edge appears with some probability
independently of all other edges. We chose $p_{edge} = 0.05$, so that
at each time $k \in K$ we had the graph $G(k) = (V, E(k))$, where
$\mathbb{P}[(i, j) \in E(k)] = 0.05$ 
for all $i$ and $j$ in each other's essential neighborhoods. 
The communication graph in this case was undirected so that
$(i, j) \in E(k)$ means that agent $i$ sends its state to agent $j$
at time $k$, and vice versa. All transmissions are received instantaneously. 
The times at which the agents sent their states to the cloud
were chosen to be randomly generated times between $k_t$ and $k_{t+1}$
which were uniformly distributed and independent of all
communications and computations. 

\begin{figure}
\centering
\includegraphics[width=2.8in]{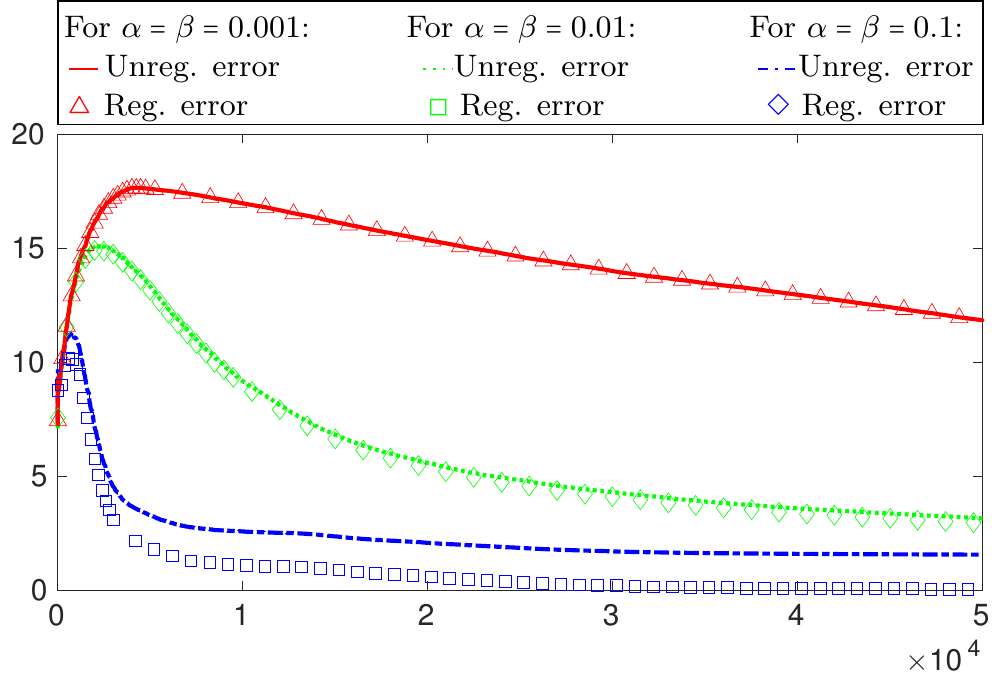}
\caption{The values of $\|x^c_t - \hat{x}\|$ (lines) and 
$\|x^c_t - \hxve\|$ (shapes) for simulation runs
using $\preg = \dreg = 0.001$ (top pair of curves),
$\preg = \dreg = 0.01$ (middle pair of curves),
and $\preg = \dreg = 0.1$ (bottom pair of curves). 
It is evident that larger regularization parameters
lead to faster decreases in error, indicating faster
convergence. 
 }
\label{fig:perror}
\end{figure}

\begin{figure}
\centering
\includegraphics[width=2.8in]{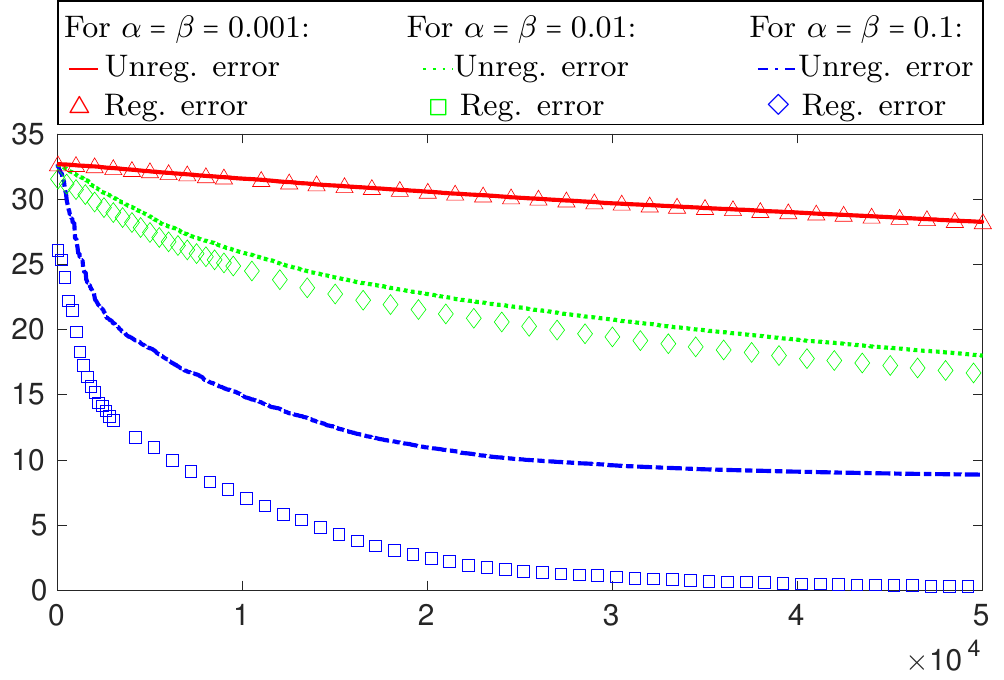}
\caption{The values of $\|\mu(t) - \hat{\mu}\|$ (lines) and 
$\|\mu(t) - \hmve\|$ (shapes) for simulation runs
using $\preg = \dreg = 0.001$ (top pair of curves),
$\preg = \dreg = 0.01$ (middle pair of curves),
and $\preg = \dreg = 0.1$ (bottom pair of curves). 
As in Figure~\ref{fig:perror}, we see that increasing
the regularization parameters $\preg$ and $\dreg$ results in
faster convergence to a final value. 
 }
\label{fig:derror}
\end{figure}

Each of the three simulation runs was run until it converged. 
In Figure~\ref{fig:perror} we see three pairs
of curves: the uppermost pair corresponds to $\preg = \dreg = 0.001$,
the middle pair corresponds to $\preg = \dreg = 0.01$, and the
lowest pair corresponds to $\preg = \dreg = 0.1$. Each pair
plots the unregularized primal error $\|x^c_t - \hat{x}\|$
for each run using lines, and the regularized primal 
error $\|x^c_t - \hxve\|$ is plotted using shapes. 
Figure~\ref{fig:derror} similarly shows the regularized and unregularized
dual errors, $\|\mu(t) - \hat{\mu}\|$ and $\|\mu(t) - \hmve\|$,
using lines and shapes, respectively, for each choice of
regularization parameters. 

\begin{table}[!tp]
\centering
\begin{tabular}{|c|c|c|c|} \hline
Value of                            & Final reg. error      & Final unreg. error    & Max final value \\ 
$\preg$ and $\dreg$                 & $\|x^c_t - \hxve\|$   & $\|x^c_t - \hat{x}\|$ & of $g_j$        \\ \hline\hline
$0.1$                               & $1.352\cdot 10^{-12}$ & $8.616$               & $1.948$         \\ \hline
$0.01$                              & $7.129\cdot 10^{-13}$ & $0.223$               & $0.252$         \\ \hline
$0.001$                             & $1.414\cdot 10^{-11}$ & $0.0237$              & $0.0262$        \\ \hline
\end{tabular} 
\caption{The final primal errors in each simulation. As predicted by Remark~\ref{rem:tradeoff},
smaller regularization parameters do indeed lead to smaller errors.}
\label{tab:perror}
\end{table}

Figure~\ref{fig:perror} shows that all error curves initially
increase, following which they decrease at different rates,
with larger regularization parameters clearly leading to faster decreases
in error. 
The final primal errors for each
simulation run are given in Table~\ref{tab:perror}, where we
see that all three runs numerically converge almost exactly to
$\hxve$. We also
see that smaller regularization parameters decrease final primal errors, as predicted by Remark~\ref{rem:tradeoff}. 

Figure~\ref{fig:derror} shows behavior in the dual space similar to
that shown in Figure~\ref{fig:perror}.
All curves appear to decrease monotonically, with larger values
of $\preg$ and $\dreg$ clearly showing a faster rate of decrease.
And as with the primal space, one finds that larger regularization
parameters lead to larger errors in the dual space; final dual error
values are shown in Table~\ref{tab:derror} wherein one finds that
all three runs virtually exactly reach $\hmve$ and, 
indeed, decreasing $\preg$ and $\dreg$ decreases the final dual error. 

\begin{table}[!tp]
\centering
\begin{tabular}{|c|c|c|} \hline
Value of                            & Final reg. error      & Final unreg. error       \\ 
$\preg$ and $\dreg$                 & $\|\mu(t) - \hmve\|$  & $\|\mu(t) - \hat{\mu}\|$ \\ \hline\hline
$0.1$                               & $7.507\cdot 10^{-12}$ & $8.616$                  \\ \hline
$0.01$                              & $4.600\cdot 10^{-12}$ & $1.573$                  \\ \hline
$0.001$                             & $1.056\cdot 10^{-10}$ & $0.174$                  \\ \hline
\end{tabular} 
\caption{The final dual errors in each simulation, which show that increasing regularization parameters
does indeed result in larger errors.}
\label{tab:derror}
\end{table}

We see in Figures~\ref{fig:perror} and \ref{fig:derror} that increasing the regularization
parameters leads to faster convergence, and this same phenomenon was observed
numerically in \cite{koshal11}. However, a key numerical difference between
our results and some of those in earlier works, e.g., \cite{boyd11} and \cite{zhang14}, is the initial
increase in distance to the optimum seen in Figure~\ref{fig:perror}. This increase is 
unavoidable due to the agents
sharing information asynchronously and is typical in simulation runs of
Algorithm $2$.

%% file: 8-conclusion.tex
\section{Conclusion} \label{sec:con}
An asynchronous multi-agent optimization algorithm for constrained
problems was presented. It was shown that the dual variable
must be kept synchronized across the agents, though their primal
updates can occur independently and with arbitrary timing.
The method presented used a Tikhonov regularization and a 
multi-agent gradient projection method to 
approximately find saddle points of the regularized Lagrangian
asynchronously. 
